\documentclass[11pt]{amsart} 
\usepackage{amsmath , amsthm , amsfonts , ifpdf}
\usepackage{color}
\usepackage{graphicx}
\usepackage{multicol}
\usepackage{epstopdf}
\usepackage{verbatim}
\usepackage[utf8x]{inputenc}

\usepackage{tikz}

\usepackage{amsmath}
\usepackage{amssymb}

\usetikzlibrary{arrows , shapes , trees , backgrounds}
\usetikzlibrary{intersections}

\theoremstyle{plain}
\newtheorem{theorem}{Theorem}[section]
\newtheorem*{theorem*}{Theorem}

\newtheorem{pro}[theorem]{Proposition}
\newtheorem{Def}[theorem]{Definition}
\newtheorem{lem}[theorem]{Lemma}
\newtheorem{Que}[theorem]{Question}

\theoremstyle{definition}
\newtheorem*{Def*}{Definition}

\newtheorem{Rem}[theorem]{Remark}

\newtheorem{exm}[theorem]{Example}

\numberwithin{equation}{section}
\newcommand{\bpo}{\begin{pro}}
\newcommand{\epo}{\end{pro}}
\newcommand{\be}{\begin{equation}}
\newcommand{\ene}{\end{equation}}
\newcommand{\br}{\begin{Rem}}
\newcommand{\er}{\end{Rem}}
\newcommand{\bl}{\begin{lem}}
\newcommand{\el}{\end{lem}}
\newcommand{\bd}{\begin{Def}}
\newcommand{\ed}{\end{Def}}
\newcommand{\ben}{\begin{enumerate}}
\newcommand{\een}{\end{enumerate}}
\newcommand{\bp}{\begin{proof}}
\newcommand{\ep}{\end{proof}}
\newcommand{\beq}{\begin{equation*}}
\newcommand{\eeq}{\end{equation*}}
\newcommand{\bear}{\begin{eqnarray*}}
\newcommand{\eear}{\end{eqnarray*}}
\newcommand{\bt}{\begin{theorem}}
\newcommand{\et}{\end{theorem}}
\newcommand{\bst}{\begin{split}}
\newcommand{\est}{\end{split}}

\newcommand{\bal}{\begin{aligned}}
\newcommand{\eal}{\end{aligned}}

\renewcommand{\P}{\partial}
\newcommand{\F}[2]{\frac{#1}{#2}}
\newcommand{\la}{\langle}
\newcommand{\ra}{\rangle}

\newcommand{\R}{\mathbb{R}}

\newcommand{\rH}{\mathrm{H}}

\newcommand{\nb}{\nabla}

\newcommand{\Mp}{maximum principle}
\newcommand{\Ode}{ordinary differential equation}

\newcommand{\wrt}{with respect to}
\newcommand{\Sc}{\varepsilon}

\newcommand{\vp}{\varphi}

\newcommand{\PLH}{{\mkern-1mu\times\mkern-1mu}}

\newcommand{\DP}{Dirichlet problem}

\newcommand{\Ca}{Caccioppoli}

\newcommand{\ap}{a priori}
\newcommand{\Eu}{Euclidean}

\newcommand{\m}[1]{\mathcal{#1}}

\newcommand{\lph}{\mathcal{L}_\phi}
\newcommand{\nf}{Nc-f}
\newcommand{\pmz}{perimeter minimizer}
\newcommand{\lfp}{locally finite perimeter set}
\newcommand{\SY}{Schoen-Yau}
\newcommand{\au}{auxiliary}
\newcommand{\xpm}{X^p_m}

\newcommand{\Suw}{\omega}
\newcommand{\HSu}{\F{Du}{\Suw}}

\def\XXint#1#2#3{{\setbox0=\hbox{$#1{#2#3}{\int}$}
		\vcenter{\hbox{$#2#3$}}\kern-.5\wd0}}

\makeatletter
\def\@citestyle{\m@th\upshape\mdseries}
\def\citeform#1{{\bfseries#1}}
\def\@cite#1#2{{%
		\@citestyle[\citeform{#1}\if@tempswa, #2\fi]}}
\@ifundefined{cite }{%
	\expandafter\let\csname cite \endcsname\cite
	\edef\cite{\@nx\protect\@xp\@nx\csname cite \endcsname}%
}{}
\makeatother
\begin{document}

\title[Dirichlet problem]{A blow-up method to prescribed mean curvature graphs with fixed boundaries} 
%\title{Some applications on $C^2$ almost minimal boundaries, III embeddedness of the Plateau problems} 
\author{Hengyu Zhou}
%\footnote{ \textbf{A data availability statement} Data sharing not applicable to this article as no datasets were generated or analysed during the current study.}

\address[H. ~Z]{ College of Mathematics and Statistics, Chongqing University, Huxi Campus, Chongqing, 401331, P. R. China}
\address{Chongqing Key Laboratory of Analytic Mathematics and Applications, Chongqing University, Huxi Campus, Chongqing, 401331, P. R. China}
\email{zhouhyu@cqu.edu.cn}
\subjclass[2010]{Primary 53A10: Secondary 35A01, 35J25}
\begin{abstract} In this paper, we apply a blow-up method of Schoen and Yau in \cite{SY81} to study a large class of prescribed mean curvature (PMC) Dirichlet problems in $n(n\geq 2)$-dimensional Riemannian manifolds. In this process we establish  curvature estimates for almost minimizing PMC hypersurfaces, using an approach of Schauder estimates from Simon \cite{Sim76}. We define an Nc-f domain, where $f$ is a given function generating from the PMC equation. Combining this condition with a sufficiently mean convex assumption the blow-up method yields corresponding solutions to these PMC Dirichlet problems. Such Nc-f assumption is almost optimal by an example. An application of our result into the PMC Plateau problem is also presented. 
	\end{abstract}
\date{\today}
\maketitle

\section{Introduction} 
In this paper we will study the prescribed mean curvature (PMC) graph in Riemannian manifolds under a fixed Dirichlet data, i.e, the Dirichlet problem of some types of PMC equations from a geometrical viewpoint. \\
\indent  The first motivation comes from the connection between PMC equations and the PMC Plateau problem of surfaces in 3-manifolds by Gulliver-Spruck \cite{GS72} and Giusti \cite{Giusti-78}. Combining their work together, in 3-manifolds, the existence of PMC graphs is closely related to the existence of PMC disks (with fixed Jordan curve). See Theorem \ref{theorem:vector:fields:two} and Theorem \ref{thm:vector:field}. \\
\indent The second motivation of our study is the PMC graph of functions to solve \eqref{eq:PMC} below which appears naturally from minimal graphs in the conformally product manifold (see Remark \ref{rk:min:wp}). In fact those functions usually satisfy a general form of PMC equations as follows: 
\begin{equation}\label{eq:PMC}
-div(\HSu)+G(x,u, -\F{Du}{\Suw},\F{1}{\Suw})=0\quad \omega=\sqrt{1+|Du|^2} 
	\end{equation}  
on a $C^2$ domain $\Omega$ in $N$. Here $(x,u)$ denotes the position  and $(-\F{Du}{\Suw},\F{1}{\Suw})$ reflects the normal vector of PMC graphs in product manifolds.\\
\indent The existence of the solution to \eqref{eq:PMC} with Dirichlet data in Riemannian manifold can be explained as to find PMC graphs in product manifolds (with a fixed boundary) as a special case to solve the PMC Plateau problem with an appropriate formulation. However, this explanation is little discussed in literature. Within our knowledge, most conditions on the underlying domain for the PMC equations \eqref{eq:PMC} for example (but not limited to) \cite{JS68,Ser69, Spr07,Tsu20,CHH19, CDO16, Bergner08,GS72,Giusti-78}, are purely analytic and have little geometric correspondence. \\
\indent To seek a more geometrical condition to investigate \eqref{eq:PMC}, we turn to the existence of Jang graphs (a special type of PMC graphs) which appear in the proof of the positive mass theorem by Schoen-Yau \cite{SY81}. They uses a blow-up method. A key observation of this method is that the non-existence of apparent horizons in the underlying domain implies the global existence of a Jang graph over this domain under the natural setting in general relativity. This inspires us to propose a new assumption to include non-existence of apparent horizons as a special case. \\
\indent Suppose $M$ is an $n$-dimensional Riemannian manifold and $f$ is a $C^{1,\alpha}$ function in the tangent bundle $TM$. When $7\geq n\geq 2$, we say a $C^2$ domain $\Omega$ in $M$ has \emph{the {\nf} property} or is \emph{an {\nf} domain} if there is no $C^{3,\alpha}$ domain $E$ in the closure of $\Omega$ with mean curvature $f(.,\vec{v})$ on the whole $\P E$ or $-f(.,-\vec{v}))$ on the whole $\P E$. Here $\vec{v}$ is the outward normal vector of $E$ and the mean curvature of $\P E$ is the divergence of $\vec{v}$. \\
\indent A simple version of our main result in Theorem \ref{thm:max:A} is stated as follows. 
\bt \label{thm:max:special}
	Let $\Omega$ be a $C^2$ compact {\nf} domain in an $n$-dimensional Riemannian manifold  where $f$ is a $C^{1,\alpha}$ function on the tangent bundle of $\Omega$ and $2\leq n\leq 7$.  Suppose  the mean curvature of the boundary of $\P \Omega$ with respect to the outward normal vector, $H_{\P\Omega}$, satisfies 
	\begin{equation}\label{eq:condition:inequality}
		H_{\P \Omega}\geq \max\{f(x,\gamma), -f(x,-\gamma)\}
		\end{equation} 
	where $\gamma$ is the outward normal vector of $\Omega$ at $x$ and the Dirichlet problem takes the form of 
	\begin{equation}\label{eq:PMC:new}
	\left\{
	\begin{split}
		-div(\F{Du}{\omega})(x)&-F(x,-\F{Du}{\omega},\F{1}{\omega})+\phi(x,u,-\F{Du}{\omega},\F{1}{\omega}) \F{1}{\omega}=0 \quad \text{ on  } \Omega \\
		u&=\psi \text{ on } \P \Omega
	\end{split}\right. 
	\end{equation}
	with $\F{\P\phi}{\P u}\geq 0$ and $F(x,X,0)=f(x,X)$ for $x\in M, X\in T_x M$. Here $F,\phi$ are $C^{1,\alpha}$ with respect to its arguments and $\phi$ is either uniformly bounded or $ \F{\P\phi}{\P u}\geq\beta>0$ for a fixed constant $\beta$. \\
 \indent Then the Dirichlet problem \eqref{eq:PMC:new} admits a unique solution in $C^{3,\alpha}(\Omega)\cap C(\bar{\Omega})$ and for any $\psi\in C(\P\Omega)$. 
	\et
\indent We give some remarks on this theorem. The condition $\F{\P \phi}{\P u}\geq 0$ guarantees the uniqueness part. In Theorem \ref{theorem:optimal} for fixed any  constant function $F$ and $\phi\equiv 0$ we construct  a compact manifold in which all conditions of Theorem \ref{thm:max:special} are satisfied except the {\nf} property, no solution to the PMC equation in \eqref{eq:PMC:new} exists. This implies that the {\nf} property is almost optimal.  Applying comparison techniques, bounded domains in Euclidean spaces, Hyperbolic spaces with sufficiently mean convex boundaries are {\nf} domains by Theorem \ref{thm:Nf}. A special case of PMC graphs from the Dirichlet problem \eqref{eq:PMC:new}, $F\equiv 0, \phi$ is a constant, has be considered by the author in \cite{GZ20} and \cite{zhou19a} with a different technique from Giusti \cite{Giu84}. We refer to a continuous method by Casteras-Heinonen-Holopainen-Lira \cite[Theorem 1.4]{Spr07} and Spruck \cite[Theorem 1.1]{CHH19} in which their assumptions are almost but not completely covered by the {\nf} property (Theorem \ref{thm:Nf}). At last, the theorem above also holds when $n>7$ by adding a singular set requirement in the definition of the {\nf} property in Definition \ref{def:ncf}. See Theorem \ref{thm:max:A} for a complete version.   \\
\indent An application of Theorem \ref{thm:max:special} is that we obtain new manifolds to solve the PMC Plateau problem in 3-manifolds in the setting of Gulliver-Spruck \cite{GS72}. More study in this topic is working on by Lizhi Chen and the author in another paper\cite{CZ23}. And the {\nf} property is an important condition we will use.\\   
\indent  In \eqref{eq:PMC:new} the term $F-\F{\phi}{\omega}$ can represent many PMC functions, for examples, constant mean curvature $c$, $tr(k)-k(\F{Du}{\omega}, \F{Du}{\omega})$( the Jang equation) for a symmetric $(0,2)$-tensor $k$ (see Remark \ref{rk:min:wp}). For some earlier results on minimal surface equations in Euclidean spaces and hyperbolic spaces we refer to \cite{JS68, Lin89, GS08}. For the Dirichlet problem of more general PMC functions, we refer to Bergner \cite{Bergner08} with a strict condition on the domain $\Omega$. The Dirichlet problem of translating mean curvature equations, $F\equiv 0$ and $\phi$ is a constant, in Euclidean spaces appears in \cite{Wang11, White15, Lopez18} to study the type-II singularity of the mean curvature flow. \\
 \indent The key ingredient to show Theorem \ref{thm:max:special} is the blow-up method in \cite{SY81} which shows the existence of Jang graphs in the positive mass theorem. Their basic idea to find the solution to \eqref{eq:PMC:new} can be summarized in the following two steps: 
 \begin{enumerate} 
 	\item First consider an auxiliary PMC equation 
 \begin{equation} \label{eq:blow:t}
 	-div(\HSu)	-F(x,\F{Du}{\omega},\F{1}{\omega})+\phi(x,u,\F{Du}{\omega},\F{Du}{\omega})\F{1}{\omega}+q(t,u,\F{1}{\omega})=0
 	\end{equation} 
 where $q$ is a smooth function with its arguments satisfying $\F{\P q}{\P u}\geq \alpha(t)>0$ and $\lim_{t\rightarrow 0}q(t,u,\F{1}{\omega})=0$ uniformly when $|u|\leq K$. Under various conditions, there is a classical solution to \eqref{eq:blow:t}, written as $u_t$;
 \item Then letting $t$ approach $0$ yields the existence of the solution to \eqref{eq:PMC:new} via a series of techniques from PDE and geometric measure theory; 
 \end{enumerate}
    In this paper we choose $q(t,u,\F{1}{\omega})$ as $\F{tu}{\omega}$ instead of $tu$ in \cite[(4.1)]{SY81}. The main difference is that the latter one requires that  \eqref{eq:condition:inequality} holds strictly to deal with the boundary value of $u$ but in our case  \eqref{eq:condition:inequality} is sufficient because of Theorem \ref{lm:bet:function}. This will give more room when we consider the blow-up process in domains with piecewise $C^2$ boundary of which mean curvature in its regular part is equal to $f$ only depending its position. This setting happens frequently in the PMC Plateau problem of surfaces in 3-manifolds discussed in our future work.\\
\indent To realize the first step above, we apply a Perron method from Eichmair \cite{Eich09} to solve \eqref{eq:blow:t} (see also Ju-Liu \cite{Ju-Liu16}). In the second step above the graph of $u_t$ can be viewed as a graph with bounded mean curvature (called as a $\Lambda$-perimeter minimizer). We apply an curvature estimate in Theorem \ref{ce:thmA} which generalizes a result of Simon \cite[Theorem 1]{Sim76} on a curvature estimate of minimal boundary. Its main method is from the Schauder estimate of elliptic equations.  For an alternative curvature estimate we refer to the work of Schoen-Simon \cite{SS81} observed by Eichmair \cite[Appendix A]{Eich10}. Our curvature estimate should have independent interest.     \\
\indent In this blow-up process an interior estimate of mean curvature type equations is also established in Theorem \ref{thm:interior:estimate}. This generalizes various special cases (see \cite{Bergner08, KS88,Kore86,Kore87,DLR16, Wang98, EM16, Spr07, Mar10} etc.) for any $C^2$ function $u$ satisfying \eqref{eq:PMC} 
where $G=G(x,u, X,r) $ is a $C^{1,\alpha}$ function satisfying a monotone condition $\F{\P G}{\P u}\geq 0$. \\
\begin{comment}Our geometric PDE method will be used to study the existence of hypersurfaces with prescribed mean curvature with (without) fixed boundaries in conformally product manifolds. 
\end{comment} 
\indent This paper is organized as follows.  In section 2 we collect some preliminary facts on $\Lambda$-{\pmz}s. In section 3 we establish the curvature estimates in Theorem \ref{ce:thmA}. In section 4 we consider an auxiliary {\DP} in Theorem \ref{thm:Ma}. In section 5 we discuss the {\nf} property and provide two examples. In section 6 we proceed the blow-up process and establish Theorem \ref{thm:max:A}. In section 7 we apply Theorem \ref{thm:max:A} into the PMC Plateau problem and extend the results of Gulliver-Spruck \cite{GS72} into Riemannian manifolds. In appendix A we collect some results on PMC equations in Riemannian manifolds. \\
 \indent This project is supported by NSFC no.
 11801046 and the  Natural Science Foundation of Chongqing, China no.  cstc2021jcyj-msxmX0430. 
 
 \section{Preliminary} A well-known fact is that PMC graphs with bounded mean curvature are $\Lambda$-{\pmz}s in product manifolds (Lemma \ref{lm:estimate}). In this section we collect some facts on perimeters needed for later usage. Our main references are the book of Maggi \cite{Mag12}, Giusti \cite{Giu84} and Simon \cite{Simon83}. In those references, most of the statements are described in Euclidean spaces. Their Riemannian versions can be easily obtained by little efforts. In what follows if no confusion we directly state them in Riemannian manifolds without proof.    \\
  \indent Fix $n\geq 2$. Let $N\subset \R^{n+k}$ be an $n$-dimensional Riemannian manifold with possible boundary and $TN$ be its tangent bundle. Let $dvol, div, \la, \ra$ denote the volume form, the divergence, and the inner product of $N$ respectively. Let $B_r(x)\subset N$ be an embedded ball $B_r(x)\subset N$ centering at $x$ with sufficiently small radius $r$. Let $W$ be an open set in $N$. We say a set $G\subset \subset W$ if the closure of $G$ is a compact set in $W$.  Let $E$ be a Borel set in $N$ and $\lambda_E$ denotes its characteristic function. Let $\mathrm{H}^l$ be the $l$-dimensional Hausdorff measure in $\R^{n+k}$. 
   \begin{Def}[P. 122, \cite{Mag12}]\label{def:finite:locally:set}Fix an open set $\Omega$ in $N$. We say $E$ is a {\lfp} or a {\Ca} set if for any compact set $K \subset \subset \Omega$
  \be
  	\sup_{spt(X)\subset\subset K}\{\int_\Omega\lambda_E div(X)dvol: \la X, X\ra \leq 1 \}<\infty
  \ene
  is finite. The perimeter of $E$ in $\Omega$  is given by 
  \be 
  P(E,\Omega)=	\sup_{spt(X)\subset\subset \Omega}\{\int_\Omega\lambda_E div(X)dvol: \la X, X\ra \leq 1 \}
  \ene 
 \end{Def}  
If two {\Ca} sets differ in a Lebesgue measure set their properties are unchanged. Thus we say such two {\Ca} sets are equivalent.  For example, if a {\Ca} set is not empty, then the volume of such set is positive. In fact we can always choose a unique representation in the equivalent class of a {\Ca} set with the property in \eqref{pro:st}. 
\begin{pro}[Proposition 3.1, \cite{Giu84}] \label{pro:choice}Suppose $E$ is a Borel set in $N$, there is a Borel set $\tilde{E}$ (differs with $E$ only a measure zero set) satisfying 
	\be \label{pro:st}
	0 <\mathrm{H}^n(\tilde{E}\cap B_r(x))< \mathrm{H}^{n}(B_r(x))
	\ene 
	for any $x$ in  $\P \tilde{E}$ and sufficiently small $r>0$. 
	\end{pro}
\br Throughout this paper $\{ x\in N: 	0 <\mathrm{H}^n(E\cap B_r(x))< \mathrm{H}^{n}(B_r(x))$  for some $r>0\}$ is called the measure-theoretical boundary of $E$, written as $\P E$. 
\er 
\begin{comment}
By Proposition 12.1 in \cite{Mag12}, $E$ is a {\Ca} set in an open set $W\subset M$ if and only if there is 
a $T N$-valued Radon measure $\mu_E$  such that 
\be 
\int_E div(X) dvol=\int_{W} \la X, d\mu_E\ra, 
\ene 
for any $C^1$ vector field $X$ in $TN$ with compact support in $W$. Let $|\mu_E|$ be the total variation of $\mu_E$. It is obvious that $spt(\mu_E)\subset \P  E$. 
\begin{Def}[Page. 167, \cite{Mag12}] The reduced boundary $\P^* E$ of a {\lfp} $E$ in $W$ is the set of those $x\in spt \mu_E$ such that the limit
	\be 
	\lim_{r\rightarrow 0+}\F{\mu_E(B_r(x))}{|\mu_E|(B_r(x))} \text{ exists and belongs to the unit tangent bundle of $W$}
	\ene 
\end{Def}
A useful property is that 
\be \label{use:property}
P(E,W)=H^{n-1}(\P^* E\cap W)
\ene 
for any open set $W$ in $N$. 
\end{comment}
  \begin{Def} \cite[Page 278]{Mag12}.  \label{Def:amb} Fix a nonnegative constant $\Lambda$. We say a {\Ca} set $E$ is a $\Lambda$-perimeter minimizer in an open set $W\subset N$ if for  any {\Ca} set $F$ satisfying $F\Delta E\subset  \subset W$ it holds that 
	\be \label{eq:variation:hold} 
	P(E, W)\leq P(F, W)+\Lambda vol(E\Delta F)
	\ene 
  Define $reg(\P E)$ as the set $\{x\in \P E: \P E $  is a $C^{1,\gamma}$ graph in a neighborhood of $x$ for some $\gamma\in (0,\F{1}{2})\}$ and $sing(\P E):=\P E\backslash reg(\P E)$.  
\end{Def} 
\br  If we choose $W$ as the embedded ball $B_r(p)$, a $\Lambda$-{\pmz} $E$ is a  $(\Lambda,r)$-{\pmz}.  By \cite[Theorem 26.5]{Mag12} the reduced boundary of $\Lambda$-{\pmz} is just $reg(\P E)$. Moreover in \cite{Mag12} all local properties of $(\Lambda, r)$-{\pmz}s hold for $\Lambda$-{\pmz}s. 
\er 
\bt [Theorem 17.7, \cite{Simon83}] \label{increasing:thm} Fix $\nu>0$.  Suppose $K$ is compact in the open set $W\subset N$ and the injective radius of $W$, $inj(W)$, is positive. Then there is a $d_0=\min\{ dist(\bar{W}, \P N), inj(W)\}$ and $\kappa$ depending only $\nu, K$ for any $C^2$ hypersurface $\Sigma$ with mean curvature $H_{\Sigma}$ satisfying $|H_{\Sigma}|\leq \nu$  it holds that 
\be 
e^{\kappa \sigma} \F{\rH^{n-1}(\Sigma\cap B_\sigma(p))}{\omega_{n-1}\sigma^{n-1}}
\ene 
is increasing for all $\sigma \in (0, d_0)$ for any $p$ in $K$. 
\et 
Usually, it is hard to detect the regular set in the boundary of a $\Lambda$-{\pmz}. But there is one exception as follows. 
\bl\label{lm:regularity:AM} Let $E$ be a {\Ca} set in a $C^2$ bounded domain $\Omega\subset N$ such that $\P E$ is tangent to $\P \Omega$ at $p$. If in a neighborhood of $p$, $E$ is a $\Lambda$-{\pmz} for some positive constants $\Lambda$. Then $p$ is contained in $reg(\P E)$. 
\el 
\bp By the Nash embedding theorem, we assume $N$ is an $n$-dimensional embedded submanifold in $\R^{n+k}$. Consider an isometry $T_\lambda(x)=\F{x-p}{\lambda}$. Following \cite[Lemma 3.3]{Tam82}  there is a sequence $\{\lambda_i\}_{i=1}^\infty$ with $\lim_{i=1}^\infty\lambda_i=0$ such that $T_\lambda(E)$ converges locally to  a minimal cone $C$ in $T_p N (=\R^n)$. However $E\subset \Omega$,  $p\in \P\Omega$ and $\Omega$ is $C^2$ near $p$. Then $C$ is contained in a half-space of $T_pN$. By \cite[Lemma 15.5]{Giu84}, $C$ is a half-space. Thus the Hausdorff density of $\P E$ at $p$ is $1$. By the Allard regularity theorem (\cite[Theorem 24.2]{Simon83}), $p$ belongs to $\P E$. 
\ep 
From Federer's dimension reduction process and Simons' minimizing cone,  the following regular property holds for  $\Lambda$-{\pmz}s. 
  \bt [Theorem 28.1,\cite{Mag12}] \label{thm:regularity}   Let $N$ be an $n$-dimensional Riemannian manifold. Suppose $E$ is a $\Lambda$-perimeter minimizer in an open set $W\subset N$. Then
\begin{enumerate}
	\item if $2\leq n\leq 7$, $sing(\P E)$ is empty;
	\item if $n=8$, $sing(\P E)$ is a discrete set;
	\item if $n\geq 9$, $\rH^s(sing(\P E))=0$ for any $s>n-8$. 
	\end{enumerate}
  \et 
  \bd\label{Def:local}	 We say a sequence of {\Ca} sets $\{E_i\}_{i=1}^\infty$ locally converges to a {\Ca} set $E$ in an open set $\Omega\subset N$ if for any Borel set $A\subset \subset \Omega$
  \be 
  \lim_{i\rightarrow +\infty}	 \int_A |\lambda_E-\lambda_{E_i}|dvol=0
  \ene 
  where $dvol$ is the volume form of $N$. 
  \ed  
  There are three good properties for the convergence of $\Lambda$-{\pmz}s. 
  \bl\label{lm:slice:convergence}  Suppose $\{E_j\}_{j=1}^\infty$ is a sequence of $\Lambda$-{\pmz}s in an open set $W\subset N$. Then
   \begin{enumerate}
   	\item \textit{(Proposition 21.13 and  Theorem 21.14 in \cite{Mag12})} there is  a $\Lambda$-{\pmz} $E$ such that a subsequence of  $\{E_j\}_{j=1}^\infty$, still denoted by $\{E_j\}_{j=1}^\infty$,  converges locally to $E$ in $W$; 
  \item  Fix $B_{r_0}(p)$ as an embedded ball in $W$. Then there is a  dense set $I$ in $(0,r_0)$ for any $r$ in $I$ such that 
  	\be \label{eq:det:two}
  	\lim_{j\rightarrow \infty}P(E_j, B_r(p))=P( E, B_r(p))
  	\ene
  	\item (Theorem 1 in \cite{Tam82} and Theorem 26.6 in \cite{Mag12}) Suppose $v$ is the outward normal vector of $\P  E$ at some point $x$ in $reg(\P E)$, $\{x_j: x_j\in \P E_j\}_{j=1}^\infty$ converges to $x$, then there is a positive integer $j_0$ such that there is a neighborhood of $x$, $W$, such that for any $j\geq j_0$,  $\P E_j$ is regular in $W$, $\{x_j\in W: x_j \in \P E_j\}_{j=1}^\infty$ converges to $x$ and $ \{v_j\}_{j\geq j_0}^\infty$ converges to $v$ in the tangent bundle $TN$. Here $v_j$ is the normal vector of $\P E_j$ at $x_j$. 
  	\end{enumerate}
  \el 
  \bp  We only show the conclusion (2). Since $\lambda_{E_j}$ converges to $\lambda_{E}$ a.e.,  the coarea formula implies that 
  \be \label{eq:delta}
  \lim_{i\rightarrow \infty} \rH^{n-1}(\P B_r(p)\cap (E_j\Delta E))=0
  \ene 
  for any $r$ in a dense set $I \subset (0,r_0)$. The semicontinuity of the perimeter implies that 
  \be \label{dqet:A}
  P( E, B_{r}(p))\leq \lim_{j\rightarrow \infty}\inf P (E_j, B_{r}(p))
  \ene 
  for any $r\in (0,r_0)$. Set $\tilde{E}_j=E\cap B_r(p)\cup E_j \backslash B_{r}(p)$,  by our definition 
  \be 
  \lim_{j\rightarrow \infty} \tilde{E}_j\Delta E_{j}=\emptyset,  
  \ene   
  Choose any $r'\in (r,r_0)$. By Definition \ref{Def:amb}, we have 
  \be \label{eq:delta:B}
  P(E_j, B_{r'}(p))\leq P(\tilde{E}_j, B_{r'}(p) )
  +\Lambda vol(E_j\Delta \tilde{E}_j ) 
  \ene 
  On the other hand by the trace formula \cite[Proposition 2.8]{Giu84} we have 
  \be \label{eq:delta:A}
  \begin{split}
  	P(\tilde{E}_j, B_{r'}(p) )&=P(E, B_{r}(p) )\\
  	&+P(E_j, B_{r'}(p)\backslash \bar{B}_{r}(p)) +\rH^{n-1}((E_j\Delta  E)\cap \P  B_r(p))
  \end{split}
  \ene 
  Combining \eqref{eq:delta:B} and \eqref{eq:delta:A},  we obtain 
  \be 
  \begin{split}
  P(E_j, B_{r}(p))&\leq P(E, B_{r}(p) )
  +\Lambda vol(E_j\Delta \tilde{E}_j)\\ &+\rH^{n-1}(E_j\Delta  E\cap \P  B_r(p))
  \end{split}
  \ene 
  Recall that we always assume $n\geq 2$. For any $r$ in $I$ by \eqref{eq:delta} and letting $j\rightarrow \infty$ we obtain 
  \be \label{dqet:B}
  \lim_{j\rightarrow +\infty }\sup P(E_j, B_{r}(p))\leq  P(E, B_{r}(p))
  \ene 
  Putting \eqref{dqet:A} and  \eqref{dqet:B} together we obtain the conclusion.   
  \ep
  Next, we give an example of $\Lambda$-perimeter minimizers for later use. For more examples see \cite[section 21.1]{Mag12}. For any function $f$ let $gr(f)$ denote the graph of $f$. Let $N\PLH \R$ denote the product manifold. We use the idea and concepts of Eichmair in \cite[Example A.1]{Eich09} to obtain the following result. 
   \bl \label{lm:estimate} Suppose $\Omega$ is a $C^2$ bounded domain in $N$,   $u$ is a $C^2$ function on $\Omega$  such that 
  \be \label{det:bound:mc:as}
  |div(\F{Du}{\omega})|\leq \mu \text{ on }\Omega,\quad  
  \ene 
  where $\omega=\sqrt{1+|Du|^2}$, $\mu$ is a fixed positive constant. Define its subgraph $U=\{(x,r):x\in \Omega, r<u(x)\}$. \\
  \indent Then there are two positive constants $\delta:=\delta(\Omega)$,  $\Lambda:=\Lambda(\Omega, \mu)$ such that  $U$ is a $\Lambda$-{\pmz} in $ \Omega_\delta\PLH \R\backslash (\P gr(u))$.  Here $\P gr(u)$ is the boundary of the graph of $u$ over $\Omega$, $\Omega_\delta:=\{ x\in N, d(x, \Omega)<\delta \}$ and the metric of $N\PLH \R$ is the product metric.
  \el 
  \br Notice that the boundary of $U$ is composed of two parts: $gr(u)$ in $\Omega\PLH\R$ and the set $\{(x,t):x\in \P\Omega, t<\P gr(u)\}$, which belongs to the boundary of the $C^2$ domain $\Omega\PLH\R$.  Our derivation does not work in a neighborhood of $\P gr(u)$. Here we do not require $u$ to be uniformly bounded on $\Omega$. \\
  \indent Comparing to \cite[Example A.1]{Eich09}, our $\Lambda$-{\pmz}s can be defined on some neighborhoods of  $\P\Omega\PLH\R\backslash \P gr(u)$. 
  \er 
  \bp  Let $d(x,y)$ be the distance between two points $x,y$ in $N$. Define $d(x)$ as the function $sign(x)d(x,\P\Omega)$ in a neighborhood of  $\P \Omega$ satisfying $sign(x)=1$ when $x\in \Omega$ and otherwise $sign(x)=-1$. Define a set $\Gamma_\delta:=\{ x\in N: -\delta < d(x)\leq  \delta \}$ such that for any $y\in \Gamma_\delta$ there is only one $x\in \P\Omega$ such that $d(x,y)=d(y,\P\Omega)$.  Since $\P\Omega$ is $C^{2}$, then so is $d(x)$.  Then we extend $d(x)$ into a function in $\Gamma_\delta \PLH\R$ as $d(q)=d(x)$ where $q=(x,r)$ in $\Gamma_\delta\PLH\R$.  Let $div$ be the divergence of $N\PLH\R$. Let $n$ be the dimension of $N$. We define an $n$-form on $\Gamma_\delta\PLH \R$
  \be \label{def:sigma:0}
  \sigma_0=Dd\llcorner dvol
  \ene 
  where $D$ and $vol$ are the gradient and the volume form of $N\PLH\R$ respectively.  
  It is easy to see that 
  \be\label{div:sigma:0}
  d\sigma_0=div(Dd)dvol,\quad |div(Dd)|\leq c(\delta) 
  \ene 
  in $\Gamma_\delta\PLH\R$. \\
  \indent  Let $\vec{v}$ be the upward normal vector of $gr(u)$ in $N\PLH\R$. We extend it into a unit vector field $X$ in $\Omega\PLH\R$ via vertical translation. Namely 
  \be 
  X(x, u(x)+t)=\vec{v}(x,u(x))
  \ene  
  Define 
  \be \label{def:sigma}
  \sigma_1=X\llcorner (dvol)
  \ene 
  By the definition of the divergence and \eqref{det:bound:mc:as} , for any point  $q=(x,u(x)+r)\in \Omega\PLH\R$ we have 
  \be \label{eq:sigma:mean:curvature}
  d\sigma_1(q)= -div(\F{Du}{\omega})dvol(q)
  \ene
 in $\Omega\PLH\R$. \\
 \indent Suppose $F$ is any {\Ca} set such that $F\Delta U\subset \subset \Omega_{\delta}\PLH\R \backslash(\P gr(u))$. Since  $F\cap (\Omega\PLH\R) \Delta U\subset \subset \Omega_\delta\PLH\R\backslash \P gr(u)$, by \cite[Theorem 3.9]{GZ21},  there is an $(n+1)$-integral  current $[[V]]$ in $\Omega_\delta\PLH \R$ such that 
  $$\P [[V]]=\P [[F\cap (\Omega\PLH\R)]]- \P [[U]],\quad  spt([[V]])\subset F\cap(\Omega\PLH\R)\Delta U$$ 
  Here $k$ is an integer and $V$ is a Borel set in $\Omega_\delta\PLH\R$. \\
  \indent Let $W$ be an open bounded set satisfying $F\Delta U\subset\subset W$ disjoint with $gr(\psi)$.  Let $D_{0}^{n}(W)$ be the collection of the $n$-smooth forms $X_\Sc$ with compact support in $W$ for some $\Sc>0$ satisfying 
 \begin{enumerate}
 	\item $\la X_\Sc, X_\Sc\ra\leq 1$, $V\subset sup(X)$; 
 	\item  $X_\Sc=\sigma_1$ in the intersection of  a neighborhood of $V$ and $\{(x,r): x\in\Omega,\quad  dist(x,\P \Omega)>\Sc\}$ for some $\Sc>0$ (if not empty); 
 	\item $X_\Sc=\sigma_0$ in the intersection of a neighborhood  of $V$ and $\{(x,r): dist(x, \P \Omega)<\F{\Sc}{2} \}$ (if not empty); 
 	\end{enumerate}
  Let $\mathcal{M}$ denote the mass of integral currents. By Definition
    \be\label{eq:A:B} 
  \begin{split}
  	&P(F\cap \Omega\PLH \R, W)=\mathcal{M}_{W}(\P [[U]]+\P  [[V]])\geq (\P [[U]]+ \P [[V]])(X_\Sc)
  	\end{split}
  \ene 
  for any $X_\Sc\in D_0^n(W)$. Notice that  $\P U$ is $C^2$ in $W$ because $W$ is disjoint with $\P gr(u)$, it is easy to find a sequence $\{X_\Sc\}_{\Sc>0}$ in $D^n_0(W)$ such that 
  \be \label{eq:A:C} 
  \lim_{\Sc\rightarrow 0}\P [[U]](X_\Sc)=P(U,W)
  \ene 
  For the same sequence $\{X_\Sc\}_{\Sc>0}$, by \eqref{eq:sigma:mean:curvature} and \eqref{det:bound:mc:as} we have 
  \be \label{eq:A:D}
  \begin{split}
 \lim_{\Sc\rightarrow 0} \P [[V]](X_\Sc)&= \lim_{\Sc\rightarrow 0} -\int_{V} d(X_\Sc) dvol\\
 & =-\int_{V}d\sigma_1\geq -\mu vol(F\cap (\Omega\PLH\R)\Delta U)
 \end{split}
  \ene   
 Combining \eqref{eq:A:B} with \eqref{eq:A:C} and \eqref{eq:A:D}  we obtain 
  \be \label{first:A}
  \begin{split}
  	P(F\cap (\Omega\PLH\R), W)\geq P(U,  W)-\mu vol((F\cap(\Omega\PLH\R))\Delta U)
  \end{split}
  \ene 
  \indent If $F\Delta U \subset \bar{\Omega} \PLH \R$, \eqref{first:A} directly gives the conclusion in the claim just letting $\Lambda=\mu$. Otherwise $F\cup (\Omega\PLH \R)$ is not equal to $\Omega\PLH \R$. In \eqref{eq:A:B}, \eqref{eq:A:C} and \eqref{eq:A:D}, we replace $F\cap (\Omega\PLH\R)$ with $F\cup (\Omega\PLH \R)$, $U$ with $\Omega\PLH \R$ and $X_\Sc$ with $\sigma_0$ respectively. The same derivation gives that 
  \be\label{first:B}
  \begin{split}
  	P(F\cup (\Omega\PLH\R), W)&\geq P((\Omega\PLH \R), W)-c(\delta) vol(F\cup(\Omega\PLH\R)\Delta (\Omega\PLH\R)))
  \end{split}
  \ene 
 Combining  \eqref{first:A} with \eqref{first:B} we obtain 
  \be 
  \begin{split}
  	&P(U, W)+P(\Omega\PLH\R, W) \leq P(F\cap (\Omega\PLH\R), W)+P(F\cup (\Omega\PLH\R), W)\\
  	&+\max\{c(\delta),\mu\}vol(F\Delta U) \leq P(F,W)+P(\Omega\PLH\R, W)+\Lambda vol(F\Delta U)
  \end{split}
  \ene 
Here $\Lambda$ is the constant $ \max\{c(\delta),\mu\}$ and in the second line we apply \cite[Lemma 15.1]{Giu84}. \\
\indent By the definition of $\Lambda$-{\pmz}s the proof is complete. 
  \ep 
  \section{Curvature estimates} 
  In this section, we extend Simon's curvature estimates \cite[Theorem 1]{Sim76} on minimal boundaries into the case of the boundary of PMC $\Lambda$-perimeter minimizers.  \\
 \indent  Let $N$ be an $n(\geq 2)$-dimensional Riemannian manifold with a metric $g$. Fix $\alpha\in (0,1)$ and an set $N_0\subset N$. By \cite[Page 44]{LP87} the H\"{o}lder space $C^{1,\alpha}(N_0)$ is the collection of $C^{1,\alpha}$ functions $u$ in a neighborhood of $N_0$, on which $C^{1,\alpha}$ norm 
 \be 
  ||u||_{C^{1,\alpha}(N_0)}=\sup_{x\in N_0} \{|u|(x)+|\nb u|(x)\}+\sup_{x,y\in N_0, y\neq x }\F{|\nb u(x)-\nb u(y)|}{d^\alpha(x,y)} \}\ene 
 is finite. Here the second supremum is taken over all $x\neq y$ in $N_0$ such that $y$ is in a normal coordinate of $x$ and $\nb u(y)$ is  the tensor at $x$ obtained from $\nb u$ at $y$ by parallel transport along the radial geodesic from $y$ to $x$.  \\
 \indent  Let $M$ be a $m(\geq 2 )$-dimensional manifold with a metric $\sigma$.  Let $TM$ be the tangent bundle with the natural Sasaki metric $\sigma_s$ (see \cite{Alb19,Sasaki58}) from $\sigma$. Therefore we have a $C^{1,\alpha}$ norm for $C^{1,\alpha}$ functions on $TM$. 
 \begin{Def}\label{def:C:norm} Let $W\subset M$ be any open set. Define $T_u M $ be the unit tangent bundle on $W$, i.e. $\{ (x,v):x\in W, v\in T_x M, \la v,v\ra \leq 1\}$. For a $C^{1,\alpha}$ function $f(x,v): TM\rightarrow \R $ is a $C^{1,\alpha}$ function,  we define $\m{N}_W(f)$ as 
 	\be \label{def:c:a:norm}
 	\begin{split}
 	\mathcal{N}_W(f)=||f||_{C^{1,\alpha}(T_uW)}
 	\end{split}
 	\ene 
 	\end{Def} 
\br \label{rk:form} In \eqref{def:c:a:norm} the definition of the Sasaki metric of $TM$ in general has a complicated form. But this definition is not hard to understand as follows.  First suppose $W\subset \R^{m}$, then any $C^{1,\alpha}$ function $f(x,v)$ on $T_u W$ can be written as $f(x^1,\cdots, x^m, v^1,\cdots, v^m)$. Then $\mathcal{N}_{W}(f)$ is proportional to 
\be \label{eq:local:form}
\begin{split}
\sup_{(x,v)\in T_u W} \{|f|&+\sum_{i=1}^m|\F{\P f}{\P x_i}|+|\F{\P f}{\P v^i}|\}(x,v)\\
&+\sum_{i=1}^m\{\sup_{(x,v),  (y,v)\in T_u W} \F{|\F{\P f}{\P x_i}(x,v)-\F{\P f}{\P x_i}(y,v)|}{|x-y|^\alpha}\\
&+\sup_{(x , w^1), (x,w^2)\in T_u W} \F{|\F{\P f}{\P v_i}(x,w^1)-\F{\P f}{\P v_i}(x,w^2)|}{|w^1-w^2|^\alpha}\}\\
\end{split}
\ene 
When $W$ is contained in a normal coordinate,  the Sasaki metric on $TW$ still takes the form
\be \label{expression:sigmas}
\sigma_s=\sum_{i,j,k,l=1}^m \sigma_{ij}dx^i dx^j+h_{kl}dv^k dv^l+w_{ki}dv^i dx^k
\ene 
It is not hard to verify that for any $f\in C^{1,\alpha}(TM)$, $\mathcal{N}_W(f)$ is still proportional to \eqref{eq:local:form} up to a constant depending only on \eqref{expression:sigmas}.
\er 
The main result of this section is stated as follows. 
\bt\label{ce:thmA}  Let $M$ be a $m(\geq 2)$-dimensional Riemannian manifold.  Fix any open set $W\subset \subset M$ and $K\subset W$ is a compact set. Fix a nonnegative constant  $\Lambda$, a positive constant $\nu$ and $\alpha\in (0,1)$.  Define 
$\mathcal{G}$ as the set of 
\begin{enumerate}
	\item  when $2\leq m\leq 7$ all $C^2$ $\Lambda$-{\pmz} $E$ in $W$ with the mean curvature of $\P E$ equal to a $C^{1,\alpha}$ function $f_E(x,\vec{v})$ satisfying $\mathcal{N}_W(f_E)\leq \nu$;
	\item  when $m\geq 8$ a sequence of $C^2$ $\Lambda$-{\pmz}s $\{E_i\}_{i=1}^\infty$ in $W$ with the mean curvature of $\P E_i$ equal to a $C^{1,\alpha}$ function $f_i(x,\vec{v})$ where $\mathcal{N}_W(f_i)\leq \nu$ and $\{E_i\}^\infty_{i=1}$ converges locally to a {\Ca} set $F$ in $W$ such that  $sing(\P F)$ is an empty set. 
\end{enumerate}Then there is a positive constant $\mu$ such that 
	\be \label{eq:1ft:conclusion}
	\sup_{E\in  \mathcal{G}}\{\max_{ \P E\cap K} |A|^2\}\leq \mu
	\ene 
	Here $A$ is the second fundamental form of $\P E$. 
\et 
Our proof mainly follows from Simon \cite[section 1]{Sim76}. 
\br Letting $\Lambda=0$, Simon \cite{Sim76} established the above result  for $2\leq m\leq 7$ and $m=8$ plus the condition that all sets $E$ in $\mathcal{G}$ are subgraphs. 
\er 
  In what follows we view $M$ as an embedded $m$-dimensional submanifold with possible boundary $\P M$ in the Euclidean space $\R^{m+k}$ for some $k$. 
 \subsection{Basic facts} The Allard regularity theorem is one of the fundamental results of integral varifolds. It says if at a point the density of multiplicity one integral varifolds is sufficiently close to $1$, such varifold in a neighborhood of this point is a $C^{1,\alpha}$ graph for some $\alpha\in (0,1)$. The statement of the original version  \cite[Theorem 8.19]{Allard72o} is too complicated for our applications. We record the  $C^2$ version of \cite[Theorem 24.2]{Simon83} without introducing the concept of varifolds. Let $B^s_\rho(a)$ be the ball $\{x:|x-a|<\rho\}$ in the $s$-dimensional {\Eu} space centered at $a$ with radius $\rho$. Let $\omega_s$ be the volume of $s$-dimensional unit ball in $\R^{s+1}$. 
\bt [Theorem 24.2, \cite{Simon83}] Fix $p>n>0$ and $k$.  \label{thm:allard:regularity} Suppose $\Sigma\subset \R^{n+k}$ is a $n$-dimensional $C^2$ embedded submanifold with empty boundary in the embedded {\Eu} ball $B^{n+k}_\rho(q)$, $q\in \Sigma$ satisfying that 
\be\label{eq:delta:con}
\F{\rH^n(B^{n+k}_\rho(q)\cap \Sigma)}{\omega_n\rho^{n}}\leq 1+\delta, \quad (\int_{ B^{n+k}_\rho(q)}|\vec{H}|^pd\Sigma)^{\F{1}{p}}\rho^{1-\F{n}{p}}\leq \delta 
\ene 
where $d\Sigma$ is the area form of $\Sigma$ and $ \vec{H} $ is the mean curvature vector of $\Sigma$ in $\R^{n+k}$.\\
\indent Then there are two constants $\delta=\delta(n,k,p), \gamma= \gamma(n,k,p)\in (0,1)$ such that  \eqref{eq:delta:con} implies the existence of the linear isometry $q_1\in \R^ 
{n+k}$ and $u\in C^{1,1-\F{n}{p}}(\bar{B}^n_{\gamma\rho}(0),\R^k)$ with $u(0)=0$, $\Sigma \cap B^{n+k}_{\gamma\rho}(q)=q_1(gr(u))\cap B^{n+k}_{\gamma\rho}(q)$ and 
\be \label{equality}
\rho^{-1}\sup_{B_{\gamma\rho}(0)} |u|+\sup_{B_{\gamma\rho}(0)}|Du|+\rho^{1-\F{n}{p}}\sup_{x,y\in B_{\gamma\rho}(0),x\neq y}\F{|Du(x)-D(y)|}{|x-y|^{1-\F{n}{p}}}\leq c\delta^{\F{n}{4}}
\ene 
where $c=c(n,k,p)$. 
\et

\bl \label{theta:lemma} Let $W$ be an open set in a $m(\geq 2)$-dimensional Riemannian manifold $M$ and $K\subset W$ be a compact set. Fix constants  $\Lambda,\nu,\alpha\in (0,1), r_0<dist(W,\P M)$. Let $\delta =\delta(m,k,p)$ be the constant given in Theorem \ref{thm:allard:regularity}. Let $\mathcal{G}$ be the set of $C^2$ $\Lambda$-{\pmz}s in Theorem \ref{ce:thmA}. 
\\
\indent Then there is a $\rho_1 =\rho_1(\Lambda,\nu,r_0,\delta,K)$ such that for any $E\in \mathcal{G}$, any $q\in K\cap \P E$ \eqref{eq:delta:con} holds on the embedded ball $B_{\rho}(q)$ for any $\rho\in (0,\rho_1]$.  
\el 
We mainly follow the idea of \cite[section 1]{Sim76}. 
\bp  Let $\delta =\delta(m,k,p)$ be the constant given in Theorem \ref{thm:allard:regularity}, $K\subset W$ is compact. Using a finite collection of balls covering $K$, by Definition \ref{Def:amb}, there is a universal constant $A$ such that 
\be 
\rH^{m-1}(\P E\cap B_{r_0}(p))\leq A
\ene 
for any $p\in K$ and any $E\in \mathcal{G}$. Here $dim M=m(\geq 2)$. Let $n=m-1$ and $k'=k+1$. By the definition of $\mathcal{G}$, for each $E\in \mathcal{G}$, its mean curvature $|H_{\P E}|\leq \nu$ on $W$. Therefore there is a $\rho_1<r_0$ such that 
\be 
(\int_{ B^{m+k}_\rho(q)}|\vec{H}|^pd (\P E) )^{\F{1}{p}}\rho^{1-\F{m-1}{p}}\leq \delta 
\ene 
for any $\rho\in (0,\rho_1)$, any $q\in K$ and any $E\in \mathcal{G}$. This is the second inequality in \eqref{eq:delta:con}.  \\
\indent We argue the first inequality in \eqref{eq:delta:con} by contradiction. Suppose no $\rho_0\in (0,\rho_1]$ as above such that the first inequality in \eqref{eq:delta:con} holds for some $E\in \mathcal{G}$ and some $\rho \in (0,\rho_0)$. Then there is a triple sequence $\{(p_j, E_j, \rho_j ):p_j\in K, E_j\in \mathcal{G},  \rho_j\in (0,\rho_1)\}$ such that 
\begin{align} \label{key:fact}
	&\rH^{m-1}(\P E_j\cap B_{\rho_j}(p_j))\geq \omega_{m-1}(1+\delta)\rho_j^{m-1}, \\
	& \lim_{j\rightarrow +\infty}p_j=z,   \lim_{j\rightarrow+\infty}\rho_j=0\label{key:fact:two}
\end{align}
 By Lemma \ref{lm:slice:convergence} $\{E_j\}_{j=1}^\infty$ converges locally to a $\Lambda$-{\pmz} $E$ and $z\in \P E$.\\
\indent From the monotonicity formula in Theorem \ref{increasing:thm}, we have 
\be \label{det:eq}
e^{\kappa (r-\rho_j)}\F{\rH^{m-1}(\P E_j\cap B_r(p_j))}{\omega_{m-1} r^{m-1}}\geq (1+\delta)
\ene 
for all $r\in [\rho_j, \rho_1)$. Here $\kappa$ is a constant depending on $K$ and $\nu$.  By \eqref{key:fact:two}
\be \label{eq:det:C}
\F{\rH^{m-1}(\P E_j\cap B_{r}(z))}{\omega_{m-1} (r-t_j)^{m-1}}\geq e^{\kappa(\rho_j+t_j-r)}(1+\delta)
\ene 
for all $r$ in $[\rho_j+t_j, \rho_1+t_j)$. Notice that  $\{\rho_j\}_{j=1}^\infty$ also converges to $0$. By the conclusion (2) of Lemma \ref{lm:slice:convergence} there is a $r_1<\rho_1$ and  a dense set $I \in (0,r_1)$ such that for any $r\in I $ we have 
\be \label{uv:det}
\F{\rH^{m-1}(\P E\cap B_{r}(z))}{\omega_{m-1} r^{m-1}}\geq e^{-\kappa r}(1+\delta )
\geq 1+\F{\delta}{2}
\ene 
When $2\leq  m\leq 7$, $sing(\P E)$ is empty in $W$ by Theorem \ref{thm:regularity}. When $m\geq 8$, $sing(\P E)$ is empty by the definition of $\mathcal{G}$ in Theorem \ref{ce:thmA}. No matter which case, $\P E$ is regular in $W$. Thus $\P E$ is embedded and $C^{1,\gamma}$ in $W$ for some $\gamma\in (0,\F{1}{2})$.   Therefore
\be \label{eq:det:three}
\lim_{r \rightarrow 0}\F{\rH^{m-1}(\P E\cap B_{r}(z))}{\omega_{m-1} r^{m-1}}=1
\ene 
This is a contradiction to \eqref{uv:det}. Thus we obtain the desirable conclusion. 
\ep 
\subsection{The proof of Theorem \ref{ce:thmA}}.
\bp Let  $W$ be an open set of $ M\subset \R^{m+k}$ and $K\subset W$ be compact. Let $\mathcal{G}$ be the collection of $C^2$ $ \Lambda$-{\pmz}s $E$ satisfying the condition (1) or (2) in Theorem \ref{ce:thmA}. \\
\indent By Lemma \ref{theta:lemma}, there is a $\rho=\rho(\Lambda,\nu,r_0,\delta,K)$ such that for any $E\in \mathcal{G}$, any $p\in K\cap \P E$, \eqref{eq:delta:con} holds on $B_{\rho}(p)$. Therefore there is a constant  $\gamma$ in $(0,1)$ such that \eqref{equality} holds on $B_{\gamma \rho}(p)$. Namely for any two points  $p_1,p_2$ in $B_{\gamma\rho}(p)\cap \P E$, 
 \be \label{uni:continuous}
| \vec{v}(p_1)-\vec{v}(p_2)|\leq C(n,k,\rho,\gamma)\delta^{\F{m-1}{4}}dist(p_1,p_2)
 \ene 
 where on the left-hand side we use the Euclidean distance in $\R^{m+k}$, in the right-hand side $dist$ is the intrinsic distance of $M$. Here $\delta$ is the fixed constant in \eqref{eq:delta:con}. \\
\indent Because $K\subset W$ is compact, the following fact is easily obtained. 
\bl\label{lm:harmonic:coordinate}   There is a positive radius $\rho_1 <\min\{r_0, \gamma\rho\}$ such that for any $p\in K$ on $B_{r_1}(p)$ there is a local coordinate  $\{u_1,\cdots, u_m\}$ such that $g_{ij}=\la \F{\P }{\P u_i} , \F{\P }{\P u_j}\ra$ satisfying the metric $g_{ij}(p)=\delta_{ij}$, $g_{mm}\equiv 1, g_{mi}=0$ for all $i=1,\cdots, m-1$ and 
\be \label{gradient:estimate}
C^{-1}I  \leq (g_{ij})(p)\leq C I \quad \quad
\rho_1^{1+\alpha}||g_{ij}||_{C^{1,\alpha}(B_{\rho_1}(p))}\leq C, 
\ene 
Here $i,j,k=1,\cdots,m$ for any $p$ in $B_{r_1}(p)$, $C$ is a constant only depending on $K$ and $I$ is the identity matrix. 
\el 
We denote $X_m^p$ by the unit vector field $\F{\F{\P }{\P u_m}}{{\sqrt{\la \F{\P}{\P u_m},  \F{\P}{\P u_m}\ra }}}$ in $B_{\rho_1}(p)$.  By Lemma \ref{lm:harmonic:coordinate} there is a uniformly positive constant $C_1$ such that 
\be \label{eqdkt} 
|X^p_m(p_1)-X^p_m(p_2)|\leq C_1 dist(p_1, p_2)
\ene 
for any $p\in K$ and any $p_1, p_2\in B_{\rho_1}(p)$. Here $|,|$ denotes the Euclidean distance in $\R^{m+k}$.  \\
\indent 
 For any $E\in \mathcal{G}$ and any $p\in \P E\cap K$, we take the coordinate $\{u^1,\cdots, u^m\}$ in $B_{\rho_1}(p)$ as in Lemma \ref{lm:harmonic:coordinate}. With a rotation we can require that  $X^p_m(p)$ is  the normal vector of $\P E$ at $p$. Let $\vec{v}_E$ denote the normal vector of $\P E$.
 Combining  \eqref{uni:continuous} and \eqref{eqdkt} yields that for any $p_1\in B_{\rho_1}(p)$, any $E\in \mathcal{G}$ and $p\in \P E\cap K$ 
 \be 
 |\la \vec{v}_E,X^p_m\ra(p_1) -\la \vec{v}_E,\xpm \ra(p)|= |\la \vec{v}_E,\xpm \ra(p_1) -1|\leq 2C_2dist(p_1, p)
 \ene 
 Here $C_2=\max \{C(n,k,\rho,\gamma)\delta^{\F{m-1}{4}}, C_1\}$. 
 \label{application:page:one}
 Since $C_2$ is a fixed constant independent of $E$, we choose $\rho_2<\rho_1$ sufficiently small such for any $E$ in $\mathcal{G}$, any $p\in K\cap  \P E$, any $q\in B_{\rho_2}(p)\cap \P E$, we have
 \be \label{eq:D:L}
   \la \vec{v}_E, X^p_m\ra(q) \geq \F{1}{2}
 \ene 
As a result $B_{\rho_2}(p)\cap \P E$ is written as $\{y, v(y)\}$ where $y=(u_1,\cdots, u_{m-1})$. Moreover $(0,v(0))$ represents the point $p$. By Lemma \ref{lm:harmonic:coordinate}, $X^m_p$ is orthogonal to $\{\F{\P }{\P u_1},\cdots, \F{\P }{\P u_{m-1}}\}$. Therefore  \be 
\vec{v}_E=\F{-g^{ij}v_i\F{\P}{\P  u_i}+X^p_m}{\sqrt{1+g^{ij}v_iv_j}},\quad\la \vec{v}_E,X^p_m\ra =\F{1 }{\sqrt{1+g^{ij}v_iv_j}}
\ene 
where $(g^{ij})=(g_{ij})^{-1}$ is a $(m-1)\PLH (m-1)$ matrix. Here $v_i, v_{ij}$ denote the first and second derivatives of $v$ with respect to the coordinate $\{u_1,\cdots, u_{m-1}\}$. \\
\indent By \eqref{gradient:estimate} and \eqref{eq:D:L},  there is a constant $C_3$ such that 
\be \label{key:gradient:estimate}
\sum_{i=1}^{m-1}v_i^2\leq C_3
\ene 
where $(u_1,\cdots, u_{m-1},v(u_1,\cdots,u_m))$ belongs to $\P E\cap B_{\rho_2}(p)$. Thus $\P E\cap B_{\rho_2}(p)$ contains a graph of $v$ over a domain 
\be
U_{C_4}:=\{(u^1,\cdots, u^{m-1}): u_1^2+\cdots+u_{m-1}^2\leq C_4\} \text{ in } \R^{m-1}
\ene 
where $C_4$ only depends on $C_3$. Applying the coordinate in Lemma  \ref{lm:harmonic:coordinate}, $v$ satisfies the following mean curvature type equation 
\be \label{stand:schauder}
\F{1}{\sqrt{1+g^{ij}v_iv_j}}(g^{ij}-\F{v^iv^j}{1+g^{ij}v_iv_j})v_{ij}=f(y, v,\F{-g^{ij}v_i\F{\P }{\P u_j} +\xpm}{\sqrt{1+g^{ij}v_iv_j}})
\ene 
over the domain $U_{C_4}$. Here $v^i=g^{ik}v_k$.  By Lemma  \ref{lm:harmonic:coordinate} and  \eqref{key:gradient:estimate}
we have 
\be \label{eq:estimate:g}
C^{-1}I\leq g^{ij}(y,v)\leq C I,\quad  \rho_1^{1+\alpha}||g_{ij}(y, v)||_{C^{1,\alpha}}\leq C \text{ on }  U_{C_4}
\ene 
Recall that the mean curvature of $\P E$ is $f(p,\vec{v})$ satisfying  $\mathcal{N}_W(f)\leq \nu$ (see Definition \ref{def:C:norm}).  By Definition \ref{def:C:norm}, Remark \ref{rk:form} and Lemma \ref{lm:harmonic:coordinate}, on $U_{C_4}$ we have 
\be \label{eq:f:u_i}
||f||_{C^{1,\alpha} (T_u W)}\leq  \nu
\ene
for $(x,p)\in T_u W$. Notice that the normal vector $\vec{v}=(p_1,\cdots, p_{m-1},p_m)$, $p_i=-g^{ij}v_i\F{\P }{\P u_j}\F{1}{\sqrt{1+g^{ij}v_iv_j}},i=1,\cdots, m-1$ and $p_m=\F{1}{\sqrt{1+g^{ij}v_iv_j}}X^p_m$. By \eqref{key:gradient:estimate}, \eqref{stand:schauder} is uniformly parabolic. Recall that $g_{ij}(p)=\delta_{ij}$. By \eqref{key:gradient:estimate}, \eqref{stand:schauder} is uniformly elliptic over $U_{C_4}$.  Because of $g_{ij}=0$ for $i,j=1,\cdots, m$ at $p$, with  \eqref{eq:estimate:g}, \eqref{eq:f:u_i} the classical Schauder estimate implies that 
\be 
 |A|^2(p)=\sum_{i=1}^{m-1}v_{ii}(0)\leq \mu (C_4,\nu, C)
\ene 
Notice that all constants $C_4, \nu, C$ are independent of any $p\in K$ and any $E\in \mathcal{G}$.\\
\indent Thus we conclude \eqref{eq:1ft:conclusion}. The proof is complete. 
\ep 
\section{An {\au} Dirichlet problem}
\indent  In this section we use the Perron method to study the following Dirichlet problem 
\be \label{T:au:DR}
\left\{
\begin{split}
\mathcal{L}_H(u)&=0 \quad \text{ on  } \Omega \\
 	u&=\psi \text{ on } \P \Omega 
\end{split}\right. 
\ene 
where  
\be \label{def:LH}
\mathcal{L}_H(u):= 	-div(\F{Du}{\omega})+ H_1(x,u, -\F{Du}{\omega}, \F{1}{\omega})+H_2(x,u, -\F{Du}{\omega}, \F{1}{\omega})\F{1}{\omega}
\ene 
 $\omega=\sqrt{1+|Du|^2}$, $\Omega$ is a domain in an $n(\geq 2)$ dimensional Riemannian manifold $N$, $M$ is the product manifold $N\PLH\R$ and $H_i(x, z, X, r)$ is a $C^{1,\alpha}$ function on $TM$ for $i=1,2$ with $ (x,z)\in M $ and $X+r\P_r\in T_{(x,z)}M$ with $X\in T_x N$.   \\
\indent The main result of this section is stated as follows. 
\bt \label{thm:Ma} Fix two constants $\beta>0$ and $\alpha\in (0,1)$.  Let $\Omega$ be a bounded $C^{2,\alpha}$ domain. Let $H_1(x,z,X,0)=f(x,z,X)$ for any $x\in N, z\in \R, X$ in $TN$. Suppose
\be \label{boundary:D}
\F{\P H_1}{\P z} +\F{\P H_2}{\P z }\geq \beta>0
\ene  and 
\be \label{serrin:type}
H_{\P\Omega}(x)\geq \max\{f(x,\psi(x), \gamma(x)),-f(x,\psi(x),-\gamma(x))\}  \quad x \in \P\Omega 
\ene 
Here $\gamma$ is the outward normal vector of $\P\Omega$. Then there is a unique solution to the Dirichlet problem \eqref{T:au:DR} in $C^{3,\alpha}(\Omega)\cap C(\bar{\Omega})
$ for any $\psi$ in $C(\P\Omega)$. 
\et
\br The Dirichlet problem \eqref{T:au:DR}-\eqref{boundary:D} can be used to blow up Jang equations. See Sakovich \cite[section 4]{Anna-21} for more references.
\er 
\br Bergner \cite{Bergner08} first considered the general form of \eqref{T:au:DR} in Euclidean spaces with only requiring that $\F{\P H_i}{\P z}$ is non-negative.  His smallness condition upon $H_1$, $H_2$ and $\Omega$ guarantees the {\ap} $C^0$ estimate of the solution to \eqref{T:au:DR} is uniformly bounded. In our setting the strict monotone condition \eqref{boundary:D} gives the corresponding {\ap} estimate. An improvement of the boundary condition in \cite{Bergner08} was later given by Marquardt \cite{Mar10}. 
\er 
The closed version of Theorem \ref{thm:Ma} is stated as follows.
\bt\label{thm:Ma:noboundary} If $\Omega$ is a closed $n(\geq 2)$-dimensional  manifold without boundary and \eqref{boundary:D} holds.  Then there is a unique solution to the problem $\m{L}_H(u)=0$ on $\Omega$ in $C^{3,\alpha}(\Omega)
$. 
\et
We did not consider the variational approach in \cite{Giu84} since the graph of the solution to the Dirichlet problem \eqref{T:au:DR} may not be written as the critical point of any area-type functional. The continuous method in \cite[section 11.3]{GT01}  did not work in our case because we can not keep \eqref{serrin:type} unchanged in the continuous process even if the $C^0$ estimate of the solution to \eqref{T:au:DR} is available by the maximum principle. Therefore the Perron theory used by Eichmair \cite{Eich09} may be the only choice to show Theorem \ref{thm:Ma}.
\subsection{The Maximum Principle} 
We collect two well-known maximum principles about mean curvature equations for later use. 
\bt
\label{thm:maximum:principal} Suppose $G(x, u, Du)$ is a $C^1$ function on $N\PLH\R\PLH TN$ satisfying $\P_u G(x, u, Du)\geq 0$. Suppose $u_1,u_2\in C^2(\Omega)\cap C(\bar{\Omega})$.  Define the elliptic operator 
\be 
Lu: =-div(\F{Du}{\sqrt{1+|Du|^2}})+G(x,u, Du) \text{ on }  \Omega
\ene 
If $L u_1(x)\geq L u_2(x)$ on $\Omega$ and $u_1(x)\geq u_2(x)$ on $\P  \Omega$, $u_1(x)\geq u_2(x)$ on $\Omega$. Moreover $u_1(x_0)=u_2(x_0)$ for some point $x_0\in \Omega$, $u_1(x)\equiv u_2(x)$. 
\et 
 Recall that for a hypersurface $\Sigma$ with its normal vector $\vec{v}$, its mean curvature $H_{\Sigma}$ is defined by $div(\vec{v})$. 
\bt [Proposition 4,\cite{EHLS16}] \label{thm:real:mp}Let $\Omega$ be a domain with $C^2$ boundary.  Let $G: TN\rightarrow \R$ be a $C^1$ function. Suppose  
\begin{enumerate}
	\item $\Sigma$ is a $C^2$ hypersurface in the complement of $\Omega$; 
	\item $\Sigma$ is tangent to $\P\Omega$ at a point $p$ in the interior of $\Sigma$ and $\P\Omega$;
	\item in an embedded ball $B$ containing $p$, for any $q\in B$, $H_{\Sigma}(q)\leq G(q, \vec{v}_1(q))$, $ H_{\P\Omega}(q)\geq G(q, \vec{v}_2(q))$;
\end{enumerate} 
where $\vec{v}_1$ and $\vec{v}_2$ denote the normal vector of $\Sigma$ and $\P\Omega$ respectively satisfying $\vec{v}_1(p)=\vec{v}_2(p)$ pointing into $\Omega$. \\
\indent Then $\Sigma$ coincides with $\P\Omega$ in a neighborhood of $p$. 
\et 
\br  The condition $\vec{v}_1(p)=\vec{v}_2(p)$ is essential because two unit spheres in Euclidean spaces can touch each other at a common point $p$ with $\vec{v}_1(p)=-\vec{v}_2(p)$. 
\er 
\subsection{The Perron Theory} We use the following definitions of the Perron method from Eichmair \cite[section 3]{Eich09}. In what follows we always assume  the operator $\m{L}_H$ is defined by \eqref{def:LH} satisfying \eqref{boundary:D}, $B_r(x)$ denotes the embedded ball centering at $x$ with radius $r$. 
\bl \label{lm:two:property} Suppose $\Omega$ is a bounded domain. Fix any $c>0$. For any $x$ in $\bar{\Omega}$ there is a continuous positive  function $r(x)$ such that the Dirichlet problem $\m{L}_H(u)=0$ in $B_r(x)$ with $u=\psi$ on $\P B_r$, $r\in (0, r(x))$ has a unique solution in $C^{3,\alpha}(B_r(x))\cap C(\bar{B}_r(x))$ for any $\psi\in C(\P B_r(x))$ satisfying $|\psi|\leq c$.
\el 
\bp  Fix any $x\in \bar{\Omega}$. First we assume that $\psi$ is in $ C^{2,\alpha}(\bar{B}_r(x))$. For any $s\in (0,1)$, let $v_s\in C^{3,\alpha}(B_r(x))\cap C(\bar{B}_r(x))$ be the solution to the Dirichlet problem 
\be \label{det:d:au:M}
\left\{
\begin{split}
-div(\F{Du}{\omega})+ s\{ G(x, u, -\F{Du}{\omega},\F{1}{\omega})	\}&=0 \quad \text{ on  } B_r(x)\\
	u&=s\psi  &\text{ on } \P B_r(x)
\end{split}\right. 
\ene 
 where $G(p, Y)$ be a $C^{1,\alpha}$ function on $TM$ given by 
$G(p,Y)=H_1(x,u, X,t)+H_2(x,u, X,t) t$ with $p=(x,u)$ and $Y=X+t\P_r$ for $X\in TN$.
Since $|\psi|\leq c$ and $\F{\P H_1}{\P z}+ \F{\P H_1}{\P z}\geq \beta$, by applying the maximum principle 
\be\label{eq:bound} 
\max_{\bar{B}_r(y)}|v_s|< c_0:=c_0(\mu_0,\beta, c)
\ene  Here $\mu_0:=\{x\in \Omega, X\in T\Omega, \la X,X\ra\leq 1,  |t|\leq 1,|u|\leq c, |H_1(x,\psi, X,t)| +|H_2(x,\psi, X,t)|\}$. Thus on $B_r(x)$ it holds that
\be \label{mc:setting:2}
s|H_1(x,v_s, -\F{Dv_s}{\omega},\F{1}{\omega})+H_2(x,v_s, -\F{Dv_s}{\omega},\F{1}{\omega})\F{1}{\omega}|\leq c_1=:c_1(c_0, \mu_0), 
\ene 
where $\omega=\sqrt{1+|Dv_s|^2}$ for any $s\in [0,1]$. Since $\bar{\Omega}$ is compact,  there is a continuous function $r(x)$ depending on $x$ such that for each $r\in (0,r(x))$, the mean curvature of $B_r(x)$ satisfies 
\be \label{mc:setting}
H_{\P B_r(x)}=\F{n-1}{r}+o(r)>c_1
\ene
 Therefore there is a positive constant $\mu:=\mu(c_0,\beta)$ such that  for any $s\in (0,1)$. 
\be \label{constant:mu:A}
\begin{split}
	s\mathcal{N}_W(G)\leq \mu
\end{split}
\ene  
where $W:=\{(x,u): x\in \bar{\Omega}, |u|\leq c_0\}$. 
By \eqref{mc:setting}, \eqref{mc:setting:2} and Theorem \ref{lm:bet:function}, there are two positive constants $\kappa, \nu$ (independent of $s\in [0,1]$) and $u_{s,\pm}=s\psi\pm \F{\log(1+\kappa d(x))}{\nu}$ such that 
\be \label{boundary:barrier}
\pm \mathcal{L}_{s}(u_{s,\pm})\geq 0 
\ene 
on the domain $\Gamma_{s,\pm}:=\{ x'\in B_{r}(x): \pm u_{s,\pm}(x')<c_0, d(x')<\F{1}{2\nu}\}$ and $\pm u_{\pm}=c_0$ on $\P \Gamma_{s,\pm}\backslash \P B_{r}(x)$. Applying Theorem \ref{thm:maximum:principal}, we obtain 
\be  
s\psi(x')- \F{\log(1+\kappa d(x'))}{\nu} \leq  v_s(x')\leq s\psi(x')+\F{\log(1+\kappa d(x'))}{\nu}
\ene 
for any  $x'$ in $\Gamma_{s,\pm}$. Here $d(x')=dist(x',\P B_r(x))$.  Because $v_{s}=s\psi$ on $\P B_r(x)$ this implies that  
\be \label{bound:c1:bd}
\max_{\P B_r(x)}|Dv_s|\leq c_3=:c_3(\kappa,\nu)
\ene 
 for any $s\in [0,1]$. By \eqref{boundary:D},  \eqref{constant:mu:A} and \eqref{bound:c1:bd}, the item (2) of Theorem \ref{thm:interior:estimate} implies that 
 \be  \label{bound:c1:all}
 \max_{\bar{ B}_r(y)}|Dv_s|\leq c_4=:c_4(c_3,\mu) \ene 
\indent 
Applying the classical Schauder estimate the $C^{1,\alpha}$ norm of  $v_s$ is uniformly bounded. By the continuous method in \cite[Theorem 11.4 and Section 11.3]{GT01} we obtain the existence of the solution $u\in C^{3,\alpha}(B_r(x))\cap C(\bar{B}_r(x))$ to the Dirichlet problem $\m{L}_H(u)\equiv 0$ on $B_{r}(x)$ with $u=\psi$ on $\P B_{r}(x)$ for any  $\psi \in C^{2,\alpha}(\bar{B}_r(x)))$ satisfying $|\psi|\leq c$.\\
\indent Now we consider the case that $\psi$ is in $C(\P B_r(x))$. Construct two monotone sequences $\{\psi_{i,k}\}_{k=1}^\infty$ for $i=1,2$ in $C^{2,\alpha}(\P B_r(x))$ which converges to $\psi$ increasingly  or decreasingly in the sense of $C(\P B_r(y))$. Let $\{u_{i,k}\}_{k=1}^\infty\in C^{3,\alpha}(B_r(x))\cap C(\bar{B}_r(x))$ be the solutions to the problem \eqref{det:d:au:M} for boundary data $\{\psi_{i,k}\}_{k=1}^\infty$.  Since $\F{\P \phi}{\P u}\geq 0$, by Theorem \ref{thm:maximum:principal}  $\{u_{i,k}\}_{k=1}^\infty,i=1,2$ are increasing and decreasing sequences respectively. By Theorem \ref{thm:interior:estimate} and Theorem \ref{thm:maximum:principal}, $\{u_{1,k}\}_{k=1}^\infty$ locally  converges  to a  function $u(x)$ in $C^{3,\alpha}(B_r(x))$ satisfying $\m{L}_H(u)=0$. Moreover, we have 
\be 
u_{1,k}\leq u\leq u_{2,k}\text{ on } B_r(x)
\ene 
Recall that $\{\psi_{i,k}\}_{k=1}^\infty$ for $i=1,2$ in $C^{2,\alpha}(\P B_r(x))$ converges to $\psi$ increasingly  or decreasingly in the sense of $C(\P B_r(x))$. Therefore $u$ belongs to $C(\bar{B}_r(x))$. \\
\indent This gives the existence in Lemma \ref{lm:two:property}. The uniqueness is from \eqref{boundary:D} and the maximum principle in Theorem \ref{thm:maximum:principal}. The proof is complete. \ep 
\begin{Def} [$\m{L}_H$ lift of a Perron subsolution]  \cite[Definition 3.2]{Eich09}. Let  $x$ be a point in $ \Omega$ and $r(x)$ be the continuous positive function in Lemma \ref{lm:two:property}.  Let $\underline{u}$ be a Perron subsolution (supersolution) of $\m{L}_H$ on $\Omega$ if the following property holds. 
	\begin{enumerate}
		\item [(P)]Fix any $r\in (0, r(x))$ and $B_r(x)$ be the embedded ball centered at $x$ with radius $r$ in $\Omega$. Let $u\in C(\bar{B}_r(x))\cap C^{2}(B_r(x))$ be the unique solution of $\m{L}_H(u)=0$ on $B_r(x)$ with $u=\underline{u}$ on $\P B_r(x)$. Then $u\geq \underline{u}$ ($u\leq \underline{u}$).
	\end{enumerate}  
	Moreover, the function $\hat{u}$ is defined by 
	\be\notag 
	\hat{\underline{u}}=\left\{\begin{split} 
		& u(x)\text{  on   }  B_r(x)\\
		& \underline{u}(x)\text{   on }  \Omega\backslash{B_r(x)}
	\end{split}
	\right.
	\ene
	is called as the $\m{L}_H$ lift of $\underline{u}$ with respect to $B_r(x)$. 
\end{Def}
\begin{Def}[\cite{Eich09}, Definition 3.3] \label{eq:def:Su}Let $\bar{u}$ be a Perron supersolution of $\m{L}_H$ on $\Omega$.  We denote by $\mathcal{S}_{\bar{u}}:=\{\underline{u}\in C(\bar{\Omega}): \underline{u}$ is a Perron subsolution of $\m{L}_H$ on $\bar{\Omega}$ with $\underline{u}\leq \bar{u}\}$ the class of all Perron subsolution of $\m{L}_H$ lying below $\bar{u}$. 
\end{Def}
By the maximum principle in Theorem  \ref{thm:maximum:principal}, we have the following properties on Perron subsolutions. 
\bl [Lemma 3.1, \cite{Eich09}] Let $\bar{u}$ and $\mathcal{S}_{\bar{u}}$ be given by Definition \ref{eq:def:Su}.  Then the following two basic properties of $\mathcal{S}_{\bar{u}}$ hold
\begin{enumerate}
	\item If $\underline{u},\underline{v}\in \mathcal{S}_{\bar{u}}$, then $\max\{\underline{u},\underline{v}\}\in \mathcal{S}_{\bar{u}}$;
	\item If $\underline{u}\in \mathcal{S}_{\bar{u}}$, $x\in \Omega$, $r\in (0, r(x))$ and $\hat{\underline{u}}$ is the $\lph$ lift of $\underline{u}$ with respect to $B_r(x)$, then $\hat{\underline{u}}\in \mathcal{S}_{\bar{u}}$. 
\end{enumerate}
\el 
\begin{Def} [\cite{Eich09}, Definition 3.4] Let $\bar{u}$ be a Perron supersolution of $\m{L}_H$ such that $\mathcal{S}_{\bar{u}}\neq \emptyset$.  The function $u^P(x):=\sup \{\underline{u}(x):\underline{u}\in \mathcal{S}_{\bar{u}}\}$ is called the Perron solution of $\m{L}_H$ on $\Omega$ with respect to $\bar{u}$.
\end{Def}
By the interior gradient estimate in Theorem \ref{thm:interior:estimate} we obtain the interior regularity of the Perron solution as follows. 
\bl\label{lm:atd}  Suppose $\Omega$ is a domain. Let $\bar{u}$ be a Perron supersolution on $\Omega$ such that $\mathcal{S}_{\bar{u}}\neq \emptyset$. Let $u^P(x)$ be the Perron solution of $\m{L}_H$ on $\Omega$ with respect to $\bar{u}$. Then $u^P\in C^{3,\alpha}(\Omega)$ and $\m{L}_H u^P=0$ on $\Omega$. 
\el
\br In this lemma we do not require that $\Omega$ is bounded.
\er 
\bp   Fix any $x_0\in \Omega$ and $r< \F{1}{4}r(x_0)$ such that $B_r(x_0)$ is an embedded ball and for any $y$ in $B_{r}(x_0)$,  $r(y)\geq \F{1}{2}r(x_0)>2r$.  Here we use the fact that $r(x)$ is a continuous function in Lemma \ref{lm:two:property}. By the definition of the Perron solution, there is a sequence $\{v_i(x)\}_{i=1}^\infty$ in $ \mathcal{S}_{\bar{u}}$ such that $\lim_{i\rightarrow+\infty}v_i(x_0)=u^P(x_0)$.\\
\indent  Now let $\{\hat{v}_i(x)\}_{i=1}^\infty$ be the $\m{L}_H$ lift of $v_i(x)$ with respect to  $B_r(x_0)$. Since $\hat{v}_i(x)$ is uniformly bounded on $B_{r}(x_0)$, Theorem \ref{thm:interior:estimate} implies that $\max_{x\in B_{r'}(x_0)}|D\hat{v}_i|$ is uniformly bounded for any $r'<r$. By the classical Schauder estimate so is the $C^{3,\alpha}$ norm of $\hat{v}_i$ on $B_{r'}(x_0)$. As a result, possible choose a subsequence,  $\{\hat{v}_i\}_{i=1}^\infty$ converges to $\hat{u}(x)$ in the locally $C^{3,\alpha}$ sense. Moreover $\m{L}_H \hat{u}=0$ on $B_r(x_0)$. By definition $\hat{u}\leq u^P$ on $B_r(x_0)$ and $\hat{u}(x_0)=u^P(x_0)$. \\
\indent Suppose there is a point $y_0$ in $ B_r(x_0)$ such that $\hat{u}(y_0)<u^P(y_0)$. Now let $\{v'_i\}_{i=1}^\infty$ in $\mathcal{S}_{\bar{u}} $ such that $v'_i\geq v_i$ for each $i$ $\lim_{i\rightarrow+\infty}v'_i(y_0)=u^P(y_0)$. Let $\hat{v}'_i$ be the $\m{L}_H$ lift of $\hat{v}_i$ with respect to  $B_{r(y_0)}(y_0)$. Since $r(y_0)>r$,  $x_0$ is contained in $B_{r(y_0)}(y_0)$.  As the above derivation in the previous paragraph, $\{\hat{v}'_i\}$ converges to a $\hat{u}'(x)$ in the local $C^2$ sense. Moreover $\m{L}_H\hat{u}'=0$ on $B_{r(y_0)}(y_0)$. Let $B$ be the set $B_{r(y_0)}(y_0)\cap  B_r(x_0)$. Since $\hat{v}'_i\geq \hat{v}_i$ for each $i$,   $\hat{u}'\geq \hat{u}$ on $B$. On the other hand,  $\hat{u}(x_0)=u^P(x_0)\geq \hat{u}'(x_0)$ and $\m{L}_H \hat{u}=\m{L}_H\hat{u}'\equiv 0$  on $B$. By Lemma \ref{lm:two:property},  $\hat{u}'\equiv \hat{u}$ on $B$. This is  a contradiction since $\hat{u}(y_0)<u^P(y_0)=\hat{u}'(y_0)$. As a result $u^P\equiv \hat{u}$ on $B_r(x_0)$. The proof is complete. 
\ep 
\subsection{Conclusion of Theorem \ref{thm:Ma} and Theorem \ref{thm:Ma:noboundary} }
\bp \textbf{The proof of Theorem \ref{thm:Ma}}. Suppose $\Omega$ is a bounded $C^{2,\alpha}$ domain. Without loss of generality we can assume $\psi\in C^{2,\alpha}(\bar{\Omega})$. The general case $\psi\in C(\P\Omega)$ follows from the same proof in the last part of Lemma \ref{lm:two:property} with a monotone approximation process.  \\
\indent Notice that $\beta>0$.  Define a finite positive number $\alpha_1$ as 
\be
\begin{split}
	\alpha_1: &=\F{1}{\beta}\max\{|H_1(x,\psi(x), X,r)|+|H_2(x,\psi(x), X,r)|: \\
	&x\in \bar{\Omega}, X\in T\Omega, \la X, X\ra \leq 1, |r|\leq 1\}+\max_{\bar{\Omega}}|\psi|+1
	\end{split}
\ene
Since $\F{\P H_1}{\P z}+\F{\P H_2}{\P  z}\geq \beta>0$, $\pm \m{L}_H(\pm \alpha_1)\geq 0$ on $\Omega$. As a result $\alpha_1$ and $-\alpha_1$ is the Perron supersolution and the Perron subsolution of $\m{L}_H u=0$ on $\Omega$ respectively.\\
\indent In Theorem \ref{lm:bet:function} let $c_0$ be the constant $\alpha_1+1$. By Theorem \ref{lm:bet:function}, there are two positive constants $\kappa, \nu$ and $v_{\pm}=\psi\pm \F{\log(1+\kappa d(x))}{\nu}$ such that 
\be \label{eq:fact:A}
\pm \m{L}_H (v_{\pm})\geq 0 
\ene 
on the domain $\Gamma_{\pm}:=\{ x\in \Omega: \pm v_{\pm}(x)<\alpha_1+1, d(x)<\F{1}{2\nu}\}$ and $\pm v_{\pm}=\alpha_1+1$ on $\P \Gamma_{\pm}\backslash \P\Omega$. Here $d(x)$ denotes the distance function from $\Omega$ which is well-defined because $\Omega$ is bounded. Define $\delta$ as the number 
\be 
\begin{split}
\delta:= \min\{ &\inf\{d(x,y): x\in \Gamma_{+},  v_{+}(x)\leq \alpha_1, y\in \P\Gamma_{+}\backslash \P\Omega\}, \\
&\inf\{d(x,y): x\in \Gamma_{-},  -v_{-}(x)\leq \alpha_1, y\in \P\Gamma_{-}\backslash \P\Omega\}
\}
\end{split}
\ene 
It is obvious that $\delta>0$. Now we define two functions 
\be
\pm u_\pm:=\left\{ 
\begin{split}&  \min\{ \pm  v_{\pm},  \alpha_1\} & \text{ on } \Gamma_{\pm}\\
	&\alpha_1   &\text{  on } \Omega\backslash \Gamma_{\pm}
\end{split}
\right.   
\ene 
By definition $u_+$ and $u_-$ are continuous on $\bar{\Omega}$ with $u_{\pm}=\psi$ on $\P\Omega$.\\
\indent Indeed $u_+$ and $u_-$ are the Perron subsolution and the Perron supersolution of $\mathcal{L}_H(u)$ respectively. Let $r(x)$ be the positive continuous function in Lemma \ref{lm:two:property}. For any $x\in \Omega$, define 
\be 
r_0(x)=\F{1}{2}\min\{r(x),dist(x,\P\Omega),\delta\}
\ene 
By the definition of $\delta$ for any fixed  $x\in \Omega$ and fixed $r<r_0(x)$, there are only two cases happen (a): $B_r(x)\subset \Gamma_{+}$, $u_+\leq v_+$ or (b): $ u_{+}\equiv  \alpha_{1}$ on $B_{r}(x)$. Then on $B_r(x)$ the $\m{L}_H$ lift of $u_+$, $\hat{u}_+$, is well-defined. If $B_r(x)\subset \Gamma_{+}$, $\hat{u}_+\leq v_+$ on $B_r(x)$ because $v_+$ is a Perron supersolution on $B_r(x)$. Otherwise $u_+=\alpha_1$ on $B_{r}(x)$, then $\hat{u}_+\leq \alpha_1\equiv u_+$. No matter which case we have $\hat{u}_+\leq \min\{v_+, \alpha_1\}=u_+$ on $\Gamma_{+}$. As a result, $u_+$ is a Perron supersolution on $\Omega$. As above, $u_-$ is a Perron subsolution on $\Omega$.  \\
\indent  Since $u_-\in \mathcal{S}_{u_+}$, by Lemma \ref{lm:atd} the Perron solution $u^P(x):=\sup \{v(x): x\in \Omega, v\in \mathcal{S}_{u_+} \}$ exists and $u^P \in C^{3,\alpha}(\Omega)$ with $\m{L}_H  u^P=0$ on $\Omega$. On the other hand, $u_-\leq  u^P\leq  u_+$ on $\bar{\Omega}$. Because  both of $u_+$ and $u_-$ are continuous on $\bar{\Omega}$ and equal to $\psi$ on $\P\Omega$. Thus $u^P$ belongs to $C(\bar{\Omega})\cap C^{3,\alpha}(\Omega)$ and is solution to  the Dirichlet problem in \eqref{T:au:DR}. The uniqueness is obvious because of \eqref{boundary:D} and the maximum principle in Theorem \ref{thm:maximum:principal}.   \\
\indent The proof of Theorem \ref{thm:Ma} is complete. 
\ep 
\bp \textbf{The proof of Theorem \ref{thm:Ma:noboundary}}. A key fact is that $\Omega$ is a closed Riemannian manifold without boundary.  Thus in what follows there is no derivation involving any boundary. Define a positive number 
\be 
\begin{split}
\alpha_2: &=\F{1}{\beta}\max\{|H_1(x, 0, X,r)|+|H_2(x,0, X,r)|: x\in \bar{\Omega},\\
&X\in T\Omega, \la X, X\ra \leq 1, |r|\leq 1\} +1
\end{split}
\ene 
where $\beta$ is from \eqref{boundary:D}.  A direct computation shows that $\pm\m{L}_H(\pm \alpha_2)\geq 0$. Thus $\alpha_2$, $-\alpha_2$ are the Perron supersolution and the Perron subsolution of $\m{L}_H(u)=0$ respectively. \\
\indent  The remainder is the same as the last part in the proof of Theorem \ref{thm:Ma}. We arrive the conclusion. 
\ep 
\section{The {\nf} property}
In this section we first discuss the relationship between the {\nf} property and the Ricci assumptions of a domain.  Then we construct two examples. The first one says that a {\nf} domain can contain a minimal hypersurface when $f\equiv 0$. The second one says that if the {\nf} assumption is almost optimal. That is if it is removed, there are a compact manifold such that there is no solution to  corresponding PMC equations.  
\bd[ the {\nf} property]\label{def:ncf} Fix $n\geq 2$, $\alpha\in (0,1)$ and $N$ is a $n$-dimensional Riemannian manifold. Suppose $f$ is a $C^{1,\alpha}$ function in the tangent bundle $TN$ and $\Omega$ is a $C^{2,\alpha}$ domain with embedded boundary. We say that $\Omega$ has  the {\nf} property if there is no {\Ca} set $E$ in the closure of $\Omega$ satisfying 
\begin{enumerate}
	\item for any $p\in \P E$, $E$ is a $\Lambda$-{\pmz} (see Definition \ref{Def:amb}) in a neighborhood of $p$ for some positive constant $\Lambda$;
	\item  the regular part $reg(\P E)$ (see Definition \ref{Def:amb}) is $C^{2,\alpha}$, embedded and orientable with its mean curvature equal to $f(x,\vec{v}(x))$ for any $x$ in  $reg(\P E)$ ($-f(x,-\vec{v}(x))$ for any $x$ in  $reg(\P E)$); 
\end{enumerate}
\indent Here $\vec{v}$ is the outward normal vector of $x$ in $reg(\P E)$, the mean curvature is defined as $div(\vec{v})$. 
\ed
\br In this definition we do not require that $\Omega$ is bounded. By Theorem \ref{thm:regularity}  the dimension of the Hausdorff dimension of $\P E \backslash reg(\P  E)$ is at most $n-8$. \\
\indent In the case $7\geq n\geq 2$, $\Omega$ is a {\nf} domain if and only if in its closure there is no $C^{2}$ domain $E$ with  mean curvature $f(x,\vec{v}(x))$ for all $x\in \P E$ or $-f(x,-\vec{v}(x))$ for all $x\in \P E$. Here $\vec{v}$ is the outward normal vector of $E$. 
\er 
\subsection{Examples of {\nf} domains}
\bt\label{thm:Nf} Suppose $\Omega$ is a $C^{2,\alpha}$ $n$-dimensional ($n\geq 2$) bounded domain with connected embedded boundary $\P \Omega$ in $N$. Let $\gamma(x)$ be the outward normal vector of $\P \Omega$ at $x$. 
Then $\Omega$ has the {\nf} property if 
\begin{enumerate}
	\item $f$ is a $C^{1,\alpha}$ function on $T N$ with 
	$\max_{x\in \bar{\Omega}, \vec{v} \in T_x N,\la \vec{v},\vec{v}\ra \leq 1}|f(x,\vec{\nu})|\leq \mu$ for some constant $\mu$ such that neither the mean curvature of $\P\Omega$ is equal to $f(.,\gamma(.))$ on the whole $\P\Omega$ or ( $-f(.,-\gamma(.))$ on the whole $\P \Omega$) . 
	\item one of the following holds: 
	\begin{enumerate} 
		\item  $H_{\P \Omega}\geq \mu,  \min_{x\in \bar{\Omega}, e\in T_xN, \la e, e\ra \leq 1}Ric(e,e)> -\F{\mu^2}{n-1}$; 
		\item $H_{\P \Omega} >\mu,  \min_{x\in \bar{\Omega}, e\in T_xN, \la e, e\ra \leq 1}Ric(e,e)\geq -\F{\mu^2}{n-1}$;
	\end{enumerate}
\end{enumerate}
\et  
\bp We argue the conclusion by contradiction. Suppose a non-empty  set $E$ is  a {\Ca} set in $\Omega$ such that (a) for any $p\in \P E$, $E$ is a $\Lambda_p$-{\pmz}  in a neighborhood of $p$ for some positive constant $\Lambda_p$;  (b) $reg(\P E)$ is $C^{2,\alpha}$ embedded such that $H_{reg(\P E)}(x)\equiv f(x,\vec{v})$ for any $x$ in $reg(\P E)$ (or $H_{reg(\P E)}(x)\equiv -f(x,-\vec{v})$ for any $x$ in $reg(\P E)$).  Here $\vec{v}$ is the normal outward vector of $reg(\P E)$. \\
\indent Now define $l=\inf_{x\in E}dist(\P\Omega, x)$ where $dist$ is the distance function induced by the metric of $N$. Let $l_m$ be the number $\sup_{x\in \bar{\Omega}}dist(x,\P\Omega)$. We have $l_m>l$ because $E$ is not empty and $vol(E)>0$ and $vol(\Omega\backslash E)>0$. First we assume $l>0$. \\
\indent We use an idea from Kasue \cite[Page, 120]{Kasue83}. Because $l\in (0,l_m)$, there is a $\Sc>0$ and a geodesic $\beta(s)$ for the arc-length parameter $s\in (-\Sc,l+\Sc]$ connecting $\beta(0)=y_0\in \P \Omega$ and $\beta(l)\in \P E$ such that  $dist(\P \Omega, \beta(t))=t$ for $t\in (0,l+\Sc)$. Let $\xi(y)$ be the inner normal vector of the point $y$ on $\P \Omega$. Define the map
\be 
\exp: \P \Omega\PLH  (-\Sc, l+\Sc)\rightarrow N, \text{ given by } \exp(y,t)=\exp_y(t\xi(y))
\ene  
Here $\exp_y$ denotes the exponential map at $y$. 
By the definition of $\beta(t)$, there is no conjugate point on $\{\beta(t):t\in [0, l+\Sc)\}$. Therefore there is a connected open set $V\subset \P \Omega$ containing $y_0$ and $\Sc_1\in (0, \Sc)$ such that  $exp$ is a diffeomorphism from $V\PLH (-\F{1}{2}\Sc_1, l+\F{1}{2}\Sc_1)$ into its image. \\
\indent  Let $V_t$ denote the image of  $V\PLH \{t\}$ under the map $\exp$ for each $t\in [0,l]$. Since $\P \Omega$ is $C^{2,\alpha}$, $V_t$ is $C^{2,\alpha}$ (see \cite{LN05}). Moreover for all $t\in (0,l)$, $V_t$ is disjoint with $\P E$ because $dist(V_t, \P \Omega)\leq t<l$. Let $H_{V_t}$ denote the mean curvature of $V_t$ at $\exp_y(t\xi(y))$ with respect the outward unit normal vector $-d(exp)(\P_t)$. We can view $\{ V_t\}_{t>0}$ as a flow with the velocity $\P_t$ as a function on $V\PLH [0,l+\F{\Sc}{2}]$. Choose $\vec{v}=-d(exp)(\P_t)$ and $f=-1$, $\P_t V_t=f\vec{v}$ and by \cite[Theorem 3.2 (v)]{HP96} we have 
\be \label{key:evolution:AN}
\begin{split}
	\F{\P H_{V_t}}{\P t}&=|A|^2+Ric(-d(exp)(\P_t)),-d(exp)(\P_t))\\
	&\geq \F{H_{V_t}^2}{n-1}-\F{\mu^2}{n-1} \text{\quad by the condition (2) }
\end{split}
\ene 
where $A$ is the second fundamental form of $V_t $ in $N$. From the condition (2) and $V_0\subset \P\Omega$, we always have $H_{\P\Omega}\geq \mu$. By the {\Mp} of {\Ode} we obtain that 
\be \label{at:d:q}
H_{V_l}\geq \mu 
\ene 
in a neighborhood of $\beta(l)$ . Since in neighborhood of $\beta(l)$ $E$ is a $\Lambda$-{\pmz}, by Lemma \ref{lm:regularity:AM} $\beta(l)$ belongs to $reg(\P E)$. Because $reg(\P E) $ is $C^{2,\gamma}$ with $H_{reg(\P E)}=\pm f(x,\pm\vec{v})$ near $\beta(l)$ for some $\gamma\in (0,\F{1}{2})$.  By Theorem \ref{thm:real:mp}, \eqref{at:d:q} and \eqref{key:evolution:AN}, $V_l$ coincides with $reg(\P E)$ near $\beta(l)$ and 
\be 
H_{V_t}\equiv \mu, \quad  Ric_{exp(t\xi(y))}(-d(exp)(\P_t), -d(exp)(\P_t))\equiv \F{\mu^2}{n-1}
\ene 
for any $t\in [0,l ]$ and $y \in V$. This is a contradiction to condition (2). \\
\indent Second we assume $l=0$. Then $\P E$ is tangent to $\P \Omega$ at some point $p\in \P\Omega$. By Lemma \ref{lm:regularity:AM}  in a neighborhood of $p$ $\P E$ is  contained in $reg(\P E)$. By the condition (2) and  Theorem \ref{thm:real:mp}  $\P\Omega$ is contained in $reg(\P E)$ and its mean curvature is equal to $f(x,\vec{v})$ near $p$ or $-f(x,-\vec{v})$ near $p$ where $\vec{v}$ is the outward normal vector of $\P\Omega$. Since $\P\Omega$ is connected, this means curvature conclusion holds on the whole $\P\Omega$.  However, this is a contradiction to condition (1).  No matter which case we will obtain the contradiction. \\
\indent Therefore $\Omega$ has the {\nf} property. The proof is complete. 
\ep 
The following example shows that the {\nf} domain can have nontrivial topology even when $f\equiv 0$. 
\begin{exm}\label{ex:ncf:cmh} Let $N$ be any $n(\geq 2)$-dimensional closed Riemannian manifold with a smooth metric $\sigma$. Fix $a>0$. Define a warped product manifold 
	\be 
\{ N\PLH (-a, a), \quad \phi^2(r)(\sigma+dr^2)\}
	\ene  
	where $\phi(r):(-a,a)\PLH\R$ is a positive smooth function satisfying 
	\be \label{det:A}
	\phi'(0)=0, \pm \phi'(\pm t)>0, t\in (0,a)
	\ene 
	By \cite[Lemma 3.1]{zhou19a} the mean curvature of the slice $N_t:=N\PLH \{t\}$ with respect to the direction $\F{\P }{\P t}$ is 
	\be \label{eq:mc:wp}
	H_{ N_t}=\F{n\phi'(t)}{\phi^2(t)}\quad \text{for any\quad }  t\in (-a,a)
	\ene 
	Therefore by \eqref{det:A} $H_{\P(N\PLH(-s,s))}>0$ with respect to the outward normal vector. Notice that $N_0$ is a minimal hypersurface in $N\PLH (-s, s)$. Indeed we claim that
	\bl Let $f\equiv 0$. For any $s\in (0,a)$, $N\PLH (-s,s)$ has the {\nf} property.
	\el
	\bp  Suppose $E\subset N\PLH (-s,s)$ is a {\Ca} set satisfying (a) for any $p\in \P E$, $E$ is a $\Lambda$-{\pmz} in a neighborhood of $p$ for some positive constant $\Lambda_p$, (b) $reg(\P E)$ is embedded, $C^2$ and minimal. Since $E$ is not empty, then $vol(E)>0$ and $vol( N\PLH (-s,s)\backslash E)>0$. Then  there are two constants  $a$ and $b$ in $[-s,s]$ such that $E\subset N\PLH[a,b]$, both $\P E\cap N_a$ and $\P E\cap N_b$ are not empty. By Lemma \ref{lm:regularity:AM}, both $\P E\cap N_a$ and $\P E\cap N_b$ are contained in $reg(\P E)$. By \eqref{det:A} and \eqref{eq:mc:wp}, at least one of the mean curvature of $N_a$ and $N_b$ is strictly positive with respect to the outward normal vector from $N\PLH (a,b)$. By Theorem \ref{thm:real:mp}, this positive mean curvature slice is contained in $reg(\P E)$ with vanishing mean curvature. This is a contradiction. \\
	\indent Therefore $N\PLH (-s,s)$ has the {\nf} property when $f\equiv 0$. 
	\ep  
\end{exm}
\subsection{Our condition is almost optimal}\label{optimal:subsection}
Now we construct an example to illustrate that the Nc-f property is almost optimal in Theorem \ref{thm:max:special} and Theorem \ref{thm:max:A}. If keep other conditions in Theorem \ref{thm:max:special} unchanged and remove the {\nf} property, there are examples that no solution to the PMC equation in \eqref{eq:PMC:new} exists. \\
\indent Recall that $n(n\geq 2)$-dimensional Euclidean space $\R^{n}$ with Euclidean metric $g_E$ takes the following warped product metric. 
\begin{equation}
	(\R^n, g_E)=\{S^{n-1}\times [0,\infty), r^2\sigma_{n-1}+dr^2\}
	\end{equation} 
where $\sigma_{n-1}$ is the standard canonical metric on $S^{n-1}$. Fix any positive $\beta>0$. Define a positive smooth function $h(r)$ by 
\begin{equation} 
	h(r)=\left\{\begin{split} 
		& r \quad  \text{for} \,\in [0, \F{n}{\beta}) \\
		& e^{k r}, \,\text{for}\, r\in [k,\infty) \\
		\end{split}\right.
	\end{equation} 
where $k$ is a positive constant strictly greater than $\max\{\beta, \F{n}{\beta}\}$. 
Now define the manifold $M$ by the set $S^{n-1}\times [0,k)$ equipped with the metric $ h^2(r)\sigma_{n-1}+dr^2$. Let $M_r$ be the domain $S^{n-1}\PLH [0,r)$.	By \cite[Lemma 3.1]{zhou19a} , the mean curvature of $\P M_r$ with respect to the metric $h^2(r)\sigma_{n-1}+dr^2$ and the outward normal vector is 
\begin{equation}
	   H_{\P M_r}=(n-1)\F{h'(r)}{h(r)}
	\end{equation} 
Therefore it is easily verified that the following three properties hold for $M$. 
\bl\label{lm:property} Let $f=\beta$ and $n\geq 2$. It holds that 
\begin{enumerate}
	\item $H_{\P M}=(n-1)k>\beta$;
	\item $H_{\P M_{\F{n-1}{\beta}}}=\beta$ and $M$ does not have the {\nf} property;
   \item $M_{\F{n}{\beta}}$ is a Euclidean ball in $\R^n$ with radius $\F{n}{\beta}$. 
	\end{enumerate}
\el
As a result $M$ satisfies the conditions of Theorem \ref{thm:max:A} except the {\nf} assumption.  
\bt\label{theorem:optimal} Let $M,\beta$ be defined as above and $n\geq 2$. There is no $C^2$ solution $u$ to solve the PMC equation $-div(\F{Du}{\omega})+\beta=0$ in the interior of $M$. Here $\omega=\sqrt{1+|Du|^2}$. 
\et  
 \bp We argue it by contradiction. Suppose there is a $C^2$ function $u$ satisfying $-div(\F{Du}{\omega})+\beta=0$ on $M$. Let $\omega_{n-1}$ denote the volume  of $S^{n-1}$. For any set $E$, $|E|$ denotes the volume of $E$. Applying the divergence 
 \begin{equation*}
 	 \begin{split} 
 		\beta |M_{\F{n}{\beta}}| =\int_{M_{\F{n}{\beta}}}div(\F{Du}{\omega})&=\int_{\P M_{\F{n}{\beta}}}\la \F{Du}{\omega},\vec{v}\ra dvol_{n-1}\\
 			&< |\P M_{\F{n}{\beta}}|\\
 			&=\omega_{n-1}(\F{n}{\beta})^{n-1}
 		\end{split} 
 \end{equation*}
Here we use the fact that $|\F{Du}{\omega}|<1$ on the interior of $M$. On the other hand, $M_{\F{n}{\beta}}$ is a Euclidean ball, $\beta |M_{\F{n}{\beta}}|=\omega_{n-1}(\F{n}{\beta})^{n-1}$. This is a contradiction. The proof is complete.

 \ep 

\section{The main result}
The main difficulty  to show Theorem \ref{thm:max:A} is that there is no $C^0$ {\ap} estimate for the solution to the Dirichlet problem \eqref{eq:PMC:new}. The main idea is to use the {\nf} property to prevent the  ``real" blow-up and obtain the existence of the finite solution as in Schoen-Yau \cite[Corollary 1]{SY81}. See the work of Eichmair \cite{Eich09, Eich10} for further applications and references.\\ 
\indent The main result of this paper is stated as follows.
 \bt  \label{thm:max:A} Fix $\alpha\in (0,1)$. Suppose $N$ is an $n (\geq 2)$-dimensional Riemannian manifold. Let $ F(x,X, \R)$ be a $C^{1,\alpha}(TN\PLH\R)$ function  satisfying that $F(x,X,0)=f(x,X)$ for any $(x,X)$ in $TN$ where $x\in N, X\in T_x N$. Suppose  $\Omega$ is a bounded  $n$-dimensional $C^2$ {\nf} domain in $N$ and the mean curvature of the boundary $\P\Omega$ satisfies 
 	\be  \label{bd:condition:B}
 	H_{\P\Omega}(x)\geq \max\{\pm f(x,\pm\gamma(x))\} 
 	\ene 
 	where $\gamma(x)$ is the outward normal vector of $\P \Omega$ at $x$.\\
 	\indent 
 Then the Dirichlet problem 
\be \label{DR:equation:A}
 \left\{
 \begin{split}
 	-div(\F{Du}{\omega})(x)&-F(x,-\F{Du}{\omega},\F{1}{\omega})+\phi(x,u,-\F{Du}{\omega},\F{1}{\omega}) \F{1}{\omega}=0 \quad \text{ on  } \Omega \\
 	u&=\psi \text{ on } \P \Omega
 \end{split}\right. 
 \ene 
 admits a unique solution in $C^{3,\alpha}(\Omega)\cap C(\bar{\Omega})$ for any $\psi(x)\in C(\P\Omega)$ where $\omega=\sqrt{1+|Du|^2}$. Here $\phi(x,z, Y,r):T(N\PLH\R) \rightarrow \R$ is a smooth function with $\F{\P \phi}{\P u}\geq 0$ which is either uniformly bounded or $\F{\P\phi}{\P u}\geq \beta$ for some positive constant $\beta>0$.  
\et 
We record the results of Spruck \cite{Spr07} on constant mean curvature equations and  Casteras-Heinonen-Holopainen \cite{CHH19} on weighted minimal graphs as follows.
\bt\label{thm:Ric:curvature} Suppose $\Omega$ is a  bounded $C^{2,\alpha}$ domain in an $n$-dimensional manifold  and $\mu$ is a positive constant satisfying 
\be \label{eq:A}
 \min_{x\in \bar{\Omega}, e\in T_x N, \la e, e\ra \leq 1}Ric(e,e)\geq -\F{\mu^2}{n-1},\quad  H_{\P\Omega}\geq \mu
\ene
then the Dirichlet problem \eqref{DR:equation:A} has a solution in $C^{2,\alpha}(\Omega)\cap C(\bar{\Omega})$
\begin{enumerate} 
	\item when $F\equiv \mu$, $\phi\equiv 0$ by Spruck \cite[Theorem 1.4]{Spr07} ;
	\item when $F=\la \F{Du}{\omega}, D m(x)\ra, \phi =h'(z)$ with $sup_{x\in \bar{\Omega}, z\in \R}\{|F|+|\phi|\}\leq \mu$ by Casteras-Heinonen-Holopainen  \cite[Theorem 1.2]{ CHH19}. 
	\end{enumerate}	
Here $m, h$ are $C^2$ functions on $\Omega$. 
\et 
For an earlier similar result we refer to Dajczer-Hinojosa-de Lira\cite{DHL08}.
\br  If the inequality in  \eqref{eq:A} holds strictly, by Theorem \ref{thm:Nf} $\Omega$ has the {\nf} property. Therefore Theorem \ref{thm:max:A} recovers partially Theorem \ref{thm:Ric:curvature}. 
\er 
\br \label{rk:min:wp} The equation \eqref{DR:equation:A} appears in the following setting. Let $M_f$ be the conformally product manifold \be 
(N\PLH (a,b), e^{2f}(\sigma+dr^2))
\ene 
where $f$ is a $C^\infty$ function on $N\PLH (a,b)$. Let $u$ is a $C^3$ function over a domain in $N$. By \cite[Lemma 3.1]{zhou19a}, the graph of $u$ is minimal in $M_f$ if and only if $u$ satisfies 
\be 	-div(\F{Du}{\omega})+n\la Df, -\F{Du}{\omega}\ra +n\F{\P f}{\P r}\F{1}{\omega}=0 \quad \text{ on  } \Omega
\ene 
Here $D$ denotes the  gradient of $f$ {\wrt} $x$ for fixed $u$. 
\er 
Now we are ready to prove Theorem \ref{thm:max:A}. 
\bp Suppose the solution $u$ to \eqref{DR:equation:A} exists. By the condition upon $\phi$ and the maximum principle, there is a positive constant $\beta_1$ only depending on $\beta$ and $\psi$ such that  $|\phi|\leq \beta_1$ for any $(x,u(x))$. Fix any $\psi \in C(\P\Omega)$. By the condition upon $\phi$ and Theorem \ref{thm:Ma}, for any $t\in (0,1)$, there is a unique solution 
$u_t$ in $C^{2,\alpha}(\Omega)\cap C(\bar{\Omega})$ solving 
\be \label{eq:Smid:A}
-div(\F{Du}{\omega})-F(x,-\F{Du}{\omega},\F{1}{\omega})+(\phi(x,u,-\F{Du}{\omega}, \F{1}{\omega})+t u)\F{1}{\omega}=0 \quad \text{ on  }\Omega
\ene 
with $u=\psi$ for any $\psi\in C(\P\Omega)$. Define 
\be \label{eq:Smid:B}
F_t(p, Y):=F(x,X,r)-(\phi(x,z,X,r)+tz)r
\ene for $t\in (0,1)$ where 
$ p=(x,z),\, Y=X+r\P_r $. It is obvious that $F_t$ is a $C^{1,\alpha}$ function on $T(N\PLH \R)$. Define  the subgraph of $u_t$ as 
\be 
U_t:=\{(x,r): x\in\Omega, r<u_t(x)\}
\ene 
\bl \label{lm:fact:Gt} It holds that 
\begin{enumerate}
	\item  for each $t\in (0,1)$ the mean curvature of $\P U_t$ in $\Omega\PLH \R$ is 
	\be 
	H_{\P U_t}(p)=F_t(p, \vec{v}), 	\ene 
	Here $p=(x,u_t(x))\in \Omega\PLH\R$ and $\vec{v}=\F{-Du_t+\P_r}{\sqrt{1+|Du_t|^2}}$ is the upward normal vector at $p$. 
	\item Define $\beta_2=\max_{\{x\in \bar{\Omega}, \la X, X\ra\leq 1, |r|\leq 1\}}\{|\psi(x)|, |F(x, X,r)+\phi(x, \psi(x), X, r)|\}$. For any $t\in (0,1)$
	\be \label{eq:fact:vip}
	|tu_t|\leq \beta_2,   \quad 
	|div(\F{Du_t}{\omega})|\leq 2\beta_2\quad  \text{ on }\Omega 
	\ene 
	where $\omega=\sqrt{1+|Du_t|^2}$. 
\end{enumerate}
\el 
\bp The item (1) follows from \eqref{eq:Smid:A} and \eqref{eq:Smid:B}. The item (2) follows from applying the maximum principle in \eqref{eq:Smid:A} and the assumption on $\phi$. 
\ep 
\indent \textbf{General Case:} \label{general:case} First we assume that Theorem \ref{thm:max:A} holds in the case of $\phi \equiv 0$.  Then there are two solutions $v_{\pm}$ in $C^{2,\alpha}(\Omega)\cap C(\bar{\Omega})$ to the Dirichlet problems given by 
\be 
-div(\F{Du}{\omega})-F(x,-\F{Du}{\omega},\F{1}{\omega})\mp\F{\beta_1+\beta_2}{\omega}=0 \quad \text{ on  } \Omega
\ene 
with $u=\psi \text{ on } \P \Omega$. Here we use the fact that 
\be 
G_{\pm}(x,X,r)=F(x,X,r)\pm(\beta_1+\beta_2)r,\quad G(x,X,0)=f(x,X)
\ene  
and for such $G_{\pm}$ the conditions in Theorem \ref{thm:max:A} still hold for the same $f$. By Theorem \ref{thm:maximum:principal} we obtain 
\be \label{eq:mid:estimate}
v_{-}\leq u_t \leq v_{+}\quad \text{ on }\quad \Omega
\ene   
for any $t\in (0,1)$. Thus 
\be \label{co:estimate}
\sup_{t\in (0,1)}|u_t|\leq \kappa
\ene 
 for a positive constant $\kappa$ independent of $t$. \\
 \indent  Let $W$ be the set $\{(x,z)\in M: x\in \bar{\Omega}, |z|\leq \kappa \}$. By \eqref{eq:Smid:B} it is obvious that 
\be 
\m{N}_W(F_t)\leq c 
\ene
where $c$ is a positive constant independent of $t\in (0,1)$. By the conclusion (1) in Theorem \ref{thm:interior:estimate} the uniformly $C^0$ estimate implies that  $|Du_t|$ is locally uniformly bounded in $\Omega$. By the classical Schauder estimate the $C^{3,\alpha}$ norm of $\{u_t\}_{t\in (0,1)}$ is locally uniformly bounded in $\Omega$. Letting $t \rightarrow 0$, a sequence of $\{u_t\}_{t>0}$ will converge to a $C^{3,\alpha}(\Omega)$ function $u$ locally as $t\rightarrow 0$. This means that $u\in C^{3,\alpha}(\Omega)$ solves  
\be \label{det:au:N}
-div(\F{Du}{\omega})-F(x,-\F{Du}{\omega},\F{1}{\omega})+\phi(x,u,-\F{Du}{\omega},\F{1}{\omega})\F{1}{\omega}=0 \quad \text{ on  }\Omega
\ene 
On the other hand from  \eqref{eq:mid:estimate} we have 
\be 
v_-\leq u\leq v_+ \quad \text{ on } \quad \Omega 
\ene 
Therefore $u$ belongs to $C(\bar{\Omega})$. Then $u$ is the solution to the Dirichlet problem \eqref{DR:equation:A} for fixed boundary data $\psi\in C(\P\Omega)$. The uniqueness follows from the condition that $\F{\P \phi}{\P u}\geq 0$. Thus we obtain Theorem \ref{thm:max:A} for general $\phi$ just assuming Theorem \ref{thm:max:A} holds for $\phi\equiv 0$.\\
\indent  \textbf{The special case $\phi \equiv 0$}. From now on we assume $\phi\equiv 0$. It suffices to establish \eqref{co:estimate} for $\phi\equiv 0$ under the {\nf} property of $\Omega$.  \\
\indent  Define the set $\Omega_a$ by 
\be\label{def:omega:1}
\Omega_a=\{ x\in N: d(x,\Omega), d(x,\Omega)<a\}
\ene 
for some positive constant $a$.  Recall that $gr(\psi):=\{(x,\psi(x)): x\in \P\Omega\}$. 
Since $\P gr(u_t)=gr(\psi)$ in $\P\Omega\PLH\R$, by Lemma \ref{lm:estimate} and \eqref{eq:fact:vip} one sees that 
\bl \label{lm:fact:U_t} There are two  positive constants $\delta:=\delta(\Omega)$ and $\Lambda:=\Lambda(\Omega,\mu)$  such that $U_t$ is a $\Lambda$-perimeter minimizer in $\Omega_\delta\PLH\R \backslash gr(\psi)$ for all $t\in (0,1)$. 
\el 
Notice that $\delta$ and $\Lambda$ are independent of $t\in (0,1)$.  By Lemma \ref{lm:slice:convergence}, there is a sequence $\{t_l\}_{l=1}^\infty$ converging to $0$ such that $\{U_{t_l}\}$ converges locally to $U_0$. Moreover $U_0$ is a  $\Lambda$-{\pmz} in $\Omega_\delta\PLH\R\backslash gr(\psi)$.  By Lemma 15.1 in \cite{Giu84} $U_0$ is the subgraph of a generalized function $u_0$ over $\Omega$. Here $u_0$ may take infinity value on  $\Omega$. We define two sets 
\be \label{def:pm}
\Omega_{\pm }=\{x\in \Omega: u_0(x)= \pm \infty\}
\ene 
 \bl  \label{lm:fact:omega} Let $\Lambda,\delta$ be the two constants in Lemma \ref{lm:fact:U_t}. It holds that
\begin{enumerate}
	\item $\Omega_{\pm}$ are $\Lambda$-{\pmz}s (see Definition \ref{Def:amb}) in $\Omega_{\delta}$; 
	\item $reg(\P\Omega_{\pm})$ are embedded, $C^{2}$ and orientable with 
	\be 
	 H_{reg(\P \Omega_{\pm})}(x)=\pm f(x,\pm \gamma(x))
	\ene 
	Here  $x\in reg(\P\Omega_{\pm})$ and $\gamma(x)$ is  outward normal vector of $x\in \P\Omega_{\pm}$ {\wrt} $\Omega_{\pm}$. 
\end{enumerate} 
\el 
\bp  Define  $T_s:N\PLH \R \rightarrow N\PLH\R$ as $ T_s(x,z)=(x,z-s)$ for any $(x,z)\in N\PLH\R$ and $s\in \R$. Since $T_s$ is an isometry of $N\PLH\R$, by Lemma \ref{lm:fact:U_t} $T_s(U_0)$ are $\Lambda$-{\pmz}s in $\Omega_\delta\PLH \R\backslash gr(\psi-s)$.  By \eqref{def:pm} and Lemma \ref{lm:slice:convergence} we have
\be\label{convergence:a}
\lim_{s\rightarrow +\infty} T_{ s}(U_{+})=\Omega_{+}\PLH \R, 
\ene 
Since $\psi\in C(\P\Omega)$ and $\Omega$ is bounded, $gr(\psi-s)$ converges to $-\infty$ uniformly as $s\rightarrow +\infty$. By Lemma \ref{lm:estimate} and Lemma \ref{lm:slice:convergence}  $\Omega_{+}\PLH\R$ is a $\Lambda$-{\pmz} in $\Omega_\delta\PLH\R$. \\
\indent Without confusion, we use $P$ to denote the perimeter in $N$ and $N\PLH \R$ in the following derivation. Suppose $\Omega_{+}$ is not a $\Lambda$-{\pmz} in $\Omega_\delta$.   Then there is a {\Ca} set $E$ in $\Omega_\delta$ satisfying $E\Delta \Omega_{+}\subset \subset \Omega_\delta$ and 
\be \label{eq:AT}
P(E, \Omega_\delta)\leq P (\Omega_+,\Omega_\delta)-\Lambda vol (E\Delta\Omega_{+})-\Sc 
\ene 
for some $\Sc>0$.  Now we define a {\Ca} set $E_T$ in $N\PLH \R$ as 
\be 
\begin{split}
	\lambda_{E_T}(p) &=\lambda_{E\PLH [-T, T]}(p), \quad p \in N\PLH [-T, T]\\
	\lambda_{E_T}(p) &=\lambda_{\Omega_{+}\PLH\R}(p) ,\quad  p\in N\PLH (-\infty, -T)\cup N\PLH (T, +\infty)
\end{split}
\ene 
Here $\lambda_{E_T}$ is the characteristic function of $E_T$. 
It is easy to see that $E_T\Delta ( \Omega_{+}\PLH \R)\subset \subset \Omega_\delta\PLH (-2T, 2T)$. Let $W_{2T}$ be the open set $\Omega_\delta \PLH (-2T, 2T)$. By \eqref{eq:AT} one sees that 
\be\notag
\begin{split}
	& P(E_T, W_{2T}) +\Lambda vol(E_T\Delta (\Omega_{+}\PLH\R))-P(\Omega_{+}\PLH \R,W_{2T} ) \\
	& \leq P(E_T, \Omega_\delta \PLH(-T,T))-P(\Omega_{+}\PLH \R,\Omega_\delta \PLH (-T,T))+2vol(\Omega_\delta)+2T\Lambda vol(\Omega_{+}\Delta E)\\
	&\leq  2T(P(E,\Omega_\delta)-P(\Omega_{+},\Omega_\delta))+2vol(\Omega_\delta)+2T\Lambda vol(\Omega_{+}\Delta E)\\
	&\leq -2\Sc T+2vol(\Omega_\delta)\leq 0
\end{split}
\ene  
for any $T>\F{vol(\Omega_\delta)}{2\Sc}$. This is a contradiction that $\Omega_{+}\PLH \R$ is a $\Lambda$-perimeter minimizer in $\Omega_\delta\PLH\R$. Thus our assumption is false and  $\Omega_{+}$ is also a $\Lambda$-{\pmz} in $\Omega_\delta$. We obtain the conclusion (1).\\
\indent  By the conclusion (1) $reg(\P \Omega_{+})\PLH\R =reg(\P (\Omega_+\PLH\R))$. Thus the conclusion (2) is equivalent that
\be \label{det:conclusion2}
H_{reg(\P \Omega_+)\PLH \R}(x, r)=F(x, \gamma(x),0)
\ene 
 Fix any point $(x,r)$ in $reg(\P \Omega_{\pm}\PLH\R)$.  Since $U_t$ locally converges to $U_0$, there is a $\{s_j\}_{j=1}^\infty$ converging to $+\infty$ and $t_j$ converging to $0$ such that 
 \be 
 \lim_{j\rightarrow \infty} T_{s_j}(U_{t_j})= \Omega_+\PLH\R
 \ene 
 in $W$ where $W$ is a bounded neighborhood of $(x,r)$ in $\Omega_\delta \PLH\R$. Moreover, there is a positive number $j_0 $ such that for any $j\geq j_0$ $W$ lies above on the subgraph of $\psi-s_j$ in $\P\Omega\PLH \R$. Thus for any $j\geq j_0$ $\P T_{s_j}(U_{t_j} )\cap W$ is contained in the graph of $u_{t_j}-s_j$ over $\Omega$. By \eqref{eq:Smid:A} the mean curvature of $\P  T_{s_j}(U_{t_j})$ is 
 \be \label{det:nde}
 H_{\P  T_{s_j}(U_{t_j})}=F(x, -\F{Du_{t_j}}{\omega},\F{1}{\omega})-t_j \F{(u_{t_j}-s_j)}{\omega}
 \ene  
 where $(x,u_{t_j}-s_j)\in W$ and $\omega=\sqrt{1+|Du_{t_j}|^2}$. Now write $F_j(x, X, r,z)=F(x, X, r)-t_j(z-s_j)r$. It is easy to see that 
 \be \label{eq:normal:A} 
 \mathcal{N}_{W}(F_j)\leq \nu
 \ene  
 independent of $j\geq j_0$. 
 For any compact set $K$ in $W\backslash sing(\P \Omega_+\PLH\R)$, By Theorem \ref{ce:thmA} there is a positive constant $\kappa$ such that 
 \be \label{eq:curvature:estimate}
 \max_{\P T_{s_j}(U_{t_j})\cap K}|A|^2_{\P T_{s_j}(U_{t_j})}\leq \kappa 
 \ene  
 for any $j\geq j_0$. 
 By Lemma \ref{lm:slice:convergence} the upward normal vector of $\{\P T_{s_j}(U_{t_j})\}$ should converge to the outward normal vector of $\Omega_{+}\PLH\R$ written as $\gamma(x)$ for $x\in reg(\P\Omega_+)$. Moreover for any point $(x,u_{t_j}-s_j)$ in $W$, $\lim_{j\rightarrow \infty }t_j(u_{t_j}-s_j) =0$. By \eqref{eq:curvature:estimate} and \eqref{det:nde}, we conclude \eqref{det:conclusion2}.\\
 \indent As for the case $\Omega_-$, we consider $\tilde{u}_t=-u_t$. Let $\tilde{U}_t$ be the subgraph of $\tilde{u}_t$ over $\Omega$ with the mean curvature of $\P \tilde{U}_t$ in $\Omega\PLH\R$ given by 
 \be 
 H_{\P \tilde{U}_t}(x,r)=\tilde{F}(x, -\F{D\tilde{u}_t}{\omega},\F{1}{\omega})-t\F{\tilde{u}_t}{\omega}
 \ene 
 By \eqref{eq:Smid:A}, $\tilde{F}(x, z, r)=-F(x, -z, r )$. As the case of $\Omega_+$ $\{\tilde{U}_t\}_{t>0}$ converges locally to $\tilde{U}_0$ which is a subgraph of $\tilde{u}_0=-u_0$ over $\Omega$. In this setting 
 \be 
 \bar{\Omega}_+=\{ \tilde{u}_0=+\infty \}=\Omega_-
 \ene 
 Arguing exactly as the argument in the case of $\Omega_+$, $\Omega_-$ is a $\Lambda$-{\pmz} in $\Omega_\delta$ and the mean curvature of $reg(\P\Omega_-)$ is 
 \be 
 H_{reg(\P\Omega_-)}(x)=\tilde{F}(x,\gamma(x),0)=-F(x,-\gamma(x),0)=-f(x,-\gamma(x))
 \ene 
 We conclude in the case of $\Omega_-$. We conclude Lemma \ref{lm:fact:omega}. 
 \ep 

   \bt\label{lm:bounded} By the {\nf} property of $\Omega$ the number $\sup_{\bar{\Omega},t\in (0,1)}|u_t|$ is finite.  
   \et 
   \bp  By  the {\nf} property of $\Omega$ and the definition of the perimeter in \eqref{def:finite:locally:set} $vol(\Omega_{\pm})=0$. Here $vol$ is the volume of $N$. \\
 \indent Suppose $\sup_{\bar{\Omega},t\in (0,1)}\{|u_t|\}$ is $+\infty$. Then there is a sequence  $\{t_i\}_{i\rightarrow \infty}$ converging to $0$ and $\{x_i\}_{i=1}^\infty\in \Omega$ converging to a point $z\in \bar{\Omega}$ such that  (a) :$u_{t_i}(x_i)>i$ or (b): $u_{t_i}(x_i)<-i$ holds. The proof of both cases is similar. Thus we only give the proof of case (a).  Let $U_i$ be the subgraph of $u_{t_i}-i$ over $\Omega$. By Lemma \ref{lm:fact:U_t} $U_i$ is a $\Lambda$-{\pmz} in $\Omega_\delta \PLH\R\backslash gr(\psi-i)$. Define $p_i=(x_i,0)$. Then $p_i\in U_i$. Fix $R>0$. Because $\psi \in C(\P\Omega)$ is uniformly bounded,  there is $i_0>0$ such that for all $i\geq i_0$ $\Omega_\delta \PLH (-R, R)$ is disjoint with the subgraph of $\psi-i$ in $\P\Omega\PLH\R$.  Thus for any $i\geq i_0$ $U_i$ are $\Lambda$-{\pmz}s in $\Omega_{\delta} \PLH (-R,R)$.  Notice that $p_i\in U_i$. From \cite[Theorem 21.11]{Mag12}  we have 
   \be 
   Vol(U_i\cap B_{r}(p_i))\geq C(n) r^{n+1}
   \ene 
   for sufficiently small $r<R$ independent of $i\geq i_0$. 
  Here $ n=dim N\geq 2$, $Vol$ is the volume of the product manifold $N\PLH \R$.  
   As $i\rightarrow +\infty$, $U_i$ should converge to the set $\Omega_+\PLH\R$. As a result, one sees that 
   \be 
   Vol((\Omega_+\PLH \R)\cap B_{R}((z,0)))\geq C(n)R^{n+1}
   \ene  
   This is a contradiction to $vol(\Omega_+)=0$. With a similar derivation, case (b) will give a contradiction to that $vol(\Omega_-)=0$. \\
   \indent Thus $\sup_{\bar{\Omega},t\in (0,1)}\{|u_t|\}$ is a finite number. We arrive at the conclusion.  
   \ep  
   Since $u_t$ is uniformly bounded on $\Omega$ for $t\in (0,1)$, as the proof of the general case of Theorem \ref{thm:max:A} in page \pageref{general:case} letting $t\rightarrow 0$ the sequence $\{u_t\}_{t\in (0,1)}$ converges to a $C^{3,\alpha}$ bounded function $u_0$ on $\Omega$ satisfying 
 \be \label{eq:step:A:D} 
 	-div(\F{Du}{\omega})-F(x,\F{-Du}{\omega}, \F{1}{\omega})=0 \quad \text{ on } \quad \Omega,\quad  \omega=\sqrt{1+|Du|^2} 
 \ene 
 To complete the proof of Theorem \ref{thm:max:A}, it remains to show that the limit of $u_0(x)$ as $x$ approaches any boundary point $z$ in $\P\Omega$ is $\psi(z)$. Suppose there is a point $x_0\in \P\Omega$ the previous statement does not hold. Recall that by Lemma \ref{lm:fact:U_t} $U_0$ is a $\Lambda$-{\pmz} in $\Omega_\delta\PLH \R\backslash gr(\psi)$. There is a point $(x_0, r_0)\in \P U_0$ with $r_0\neq \psi(x_0)$. Since $\psi$ is continuous on $\P\Omega_0$ and $\Omega\PLH\R$ is a $C^2$ domain, by Lemma \ref{lm:regularity:AM} $\P U_0$ is regular in $W$ which is a sufficiently small neighborhood of $(x_0,r_0)$. As in \eqref{eq:normal:A}  and \eqref{eq:curvature:estimate}, in $W$ $\P U_0$ is $C^{3,\alpha}$ with its mean curvature equal to $F(x, X, \la \vec{v},\P_r\ra)$. Here $\vec{v}$ is the outward normal vector of $\P U_0$ in $(x,r)\in W$ with $\vec{v}=X+\la \vec{v},\P_r\ra\P_r $ for $X\in T_x N$.  At $(x_0,r_0)$ there are two cases $\vec{v}=\pm \gamma(x_0)$. However by the condition (i) in Theorem \ref{thm:max:A} we have 
 \be 
 H_{\P\Omega\PLH\R}(x,z)\geq \pm F(x,\pm \gamma(x),0 )
 \ene 
 No matter which case, applying the maximum principle in Theorem \ref{thm:real:mp}, we obtain that $\P U_0$ coincides with $\P\Omega\PLH\R$ in $W$. Since we can repeat the above process to obtain that $\P U_0$ coincides with $\Omega\PLH\R$ far from $gr(\psi)$. Namely, the boundary of the graph of $u_0$ is unbounded. This is a contradiction because $\{u_t\}_{t\in (0,1)}$ is uniformly bounded. \\
\indent Therefore $u_0\in C(\bar{\Omega})$ and $\lim_{x\in \Omega, x\rightarrow z}u_0(x)=\psi(z)$ for any $z\in \P\Omega$. With \eqref{eq:step:A:D} $u_0$ belongs to $C(\bar{\Omega})\cap C^{3,\alpha}(\Omega)$ solving the Dirichlet problem \eqref{DR:equation:A} satisfying  $u_0=\psi$ on $\P \Omega$.  The uniqueness is from the maximum principle. We conclude Theorem \ref{thm:max:A} in the case that $\phi\equiv 0$ .\\
\indent  The proof of Theorem \ref{thm:max:A} is complete.
\ep   
\section{The PMC Plateau problem of disks}
In this section we discuss the application of Theorem \ref{thm:max:A} into the PMC Plateau problem in 3-manifolds. We mainly answer the following question with the {\nf} property.  
\begin{Que}\label{question:PMC:Plateau} Let $\Gamma$ be a null-homotopic Jordan curve in a $C^2$ compact 3-manifold $M$ with connected boundary $\P M$ and $f(x)$ be a $C^{1,\alpha}$ function on $M$. Can we find an immersed disk with boundary $\Gamma$ in 3-manifold? 
	\end{Que}
\subsection{Existence}
The main result of this subsection is stated as follows. 
\bt\label{thm:max:B} Let $\Gamma, f, M$ be given as in Question \ref{question:PMC:Plateau}. Suppose $M$ is an {\nf} domain with the mean curvature of its boundary satisfying 
\begin{equation}\label{eq:condition:mc}
	H_{\P M}(x)\geq  |f(x)| \quad \forall x\in \P M
\end{equation}
Then there is an immersed disk solution to the PMC Plateau problem with boundary $\Gamma$ in the sense of Definition \ref{thm:existence} below. 
\et
\br Gulliver \cite{Gul73} studied the PMC Plateau problem of disks in a star-shaped compact 3-manifold satisfying a Ricci and mean curvature conditions. The author thought that the result of Gulliver-Spruck \cite{GS72} includes that of \cite{GS72}. The result of Duzaar-Steffen \cite{DS99} on PMC Plateau problems is different with \cite{Gul73}, \cite{GS72} and Theorem \ref{thm:max:B} in the following point: to make the PMC functional in \eqref{def:functional} non-negative,  \cite{DS99} first fixes the domain and put appropriate condition on PMC functions; in our setting we reverse the sequence. 
\er 
\indent Let $f$ be the $C^{1,\alpha}$ function given as above. First we recall the definition of PMC functionals on 3-manifold $M$ from Gulliver-Spruck \cite{GS72}. First assume there is a $C^1$ vector $Q$ satisfying 
\begin{equation} \label{eq:condition:Q}
	div(Q)=f\quad \la Q, Q\ra < 1 \quad\text{ on the interior of }\quad M
	\end{equation}  Let $B$ denote the unit open disk $\{(u,v):u^2+v^2< 1\}$. Fix a null-homotopic Jordan curve $\Gamma$ in $\P M^3$. Define $\mathcal{B}_\Gamma$ be the set of functions $X:B\rightarrow M^3$ in $C(\bar{B})\cap W^{1,2}(B)$ such that $X$ maps $\P B$ monotonely into $\Gamma$. Consider the following functional 
\begin{equation}\label{def:functional}
	\mathcal{F}(X)=\F{1}{2} \int_{B} \la X_u, X_u\ra +\la X_v,X_v\ra dudv+\int_{B} \la Q, X_u\wedge X_v\ra du dv 
\end{equation}
where $X\in \mathcal{B}_\Gamma$, $X_u=X_*(\F{\P}{\P u})$ and  $X_v=X_*(\F{\P}{\P v})$ and $Q$ is a $C^1$ vector field on $M$ satisfying $div(Q)=f(x)$.
\br Because we define the mean curvature of a smooth hypersurface by $div(\vec{v})$, our choice of $Q$ differs that of \cite{GS72} by a constant. 
\er  Consider the minimizing problem 
\begin{equation}\label{minimize:problem}
	\min_{X\in \mathcal{B}_\Gamma} \mathcal{F}(X) 
\end{equation}
\begin{Def}\label{thm:existence} Suppose there is a map $X\in \mathcal{B}_\Gamma$ realizes the minimum in \eqref{minimize:problem}. Such $X$ is called the disk solution to the PMC Plateau problem in \eqref{minimize:problem}.   
\end{Def}
For more geometric properties of $X(B)$ we refer to the work of Gulliver in \cite{Gul70}. 
On the other hand, the existence of $Q$ is the key ingredient of the PMC Plateau problem. Actually the main result of Gulliver-Spruck \cite{GS72} can be summarized as the following existence. 
\bt \label{thm:vector:field} (\cite[Proposition 4.1]{GS72}) Let $\Gamma, M, f$ be given in Question \ref{question:PMC:Plateau}. Suppose there is a $C^1$ vector field $Q$ on $M$ satisfying \eqref{eq:condition:Q} and $\la Q,Q\ra <1-\beta$ for some positive constant $\beta$. Then there is a disk solution to the PMC Plateau problem with boundary $\Gamma$ in the sense of Definition \ref{thm:existence}. 
\et 

First we show that the equivalence between PMC graphs and a $C^1$ vector field $Q$ satisfying \eqref{eq:condition:Q}.  Then it holds that 
\bt\label{theorem:vector:fields:two} Fix $M, f(x)$ as in Question \ref{question:PMC:Plateau}. The existence of a $C^1$ vector $Q$ on $M$ satisfying \eqref{eq:condition:Q} is equivalent to the existence of the graph of a $C^{3,\alpha}$ function $u$ on $M$ satisfying the PMC equation $-div(\F{Du}{\omega})+f(x)=0$ on $M$ where $\omega=\sqrt{1+|Du|^2}$. 
\et 
\bp First we show the right hand-side implies the left hand-side. If the function $u$ satisfying the PMC equation above, define $Q=\F{Du}{\omega}$. Since $u$ is $C^{3,\alpha}$, then $\la Q, Q\ra <1$ and $div(Q)=f(x)$. \\
\indent Suppose there is a $C^1$ vector field $Q$ satisfying $\la Q, Q\ra<1$ on $M$ and $div (Q)=f(x)$. By the classical divergence theorem, for any $\Ca$ set $A\neq \emptyset, M$, it holds that 
\begin{equation} 
	|\int_{A}f(x)dvol|<P(A, M')	
\end{equation}
where $M'$ is an open Riemannian manifold strictly containing  $M$ and $P$ is the perimeter of $A$ in $M'$, $dvol$ is the volume form of $M$. Notice that all conclusions in \cite{Giusti-78} hold on Riemannian manifolds without changing any words. From \cite[Theorem 1.1]{Giusti-78}, there is a $C^2$ solution $u$ to the PMC equation $-div(\F{Du}{\omega})+f=0$ on the interior of $M$. The $C^{3,\alpha}$ property of $u$ follows from the classical regularity theory when $f$ is $C^{1,\alpha}$. The proof is complete. 
\ep 
Now we are ready to show Theorem \ref{thm:max:B}. 
\bp By Theorem \ref{thm:vector:field}, it is suffices to show that there is a vector field $Q$ satisfying \eqref{eq:condition:Q} and $\la Q, Q\ra <1-\beta$ on $M$. Since $M$ is {\nf} and the mean curvature of its boundary satisfies \eqref{eq:condition:mc}, given any constant $e$, by theorem \ref{thm:max:A} there is a $C^{3,\alpha}$ function $u\in C^{3,\alpha}(int(M))\cap C(M)$ solving $-div(\F{Du}{\omega})+f=0$ in the interior of $M$ and $u\equiv e$ on $\P M$. Define 
$Q=\F{Du}{\omega}$ where $\omega=\sqrt{1+|Du|^2}$.\\
\indent  If $\sup_{x\in int(M)}\la Q, Q\ra=1$. There is a sequence of $\la x_k\ra_{k=1}^\infty$ in $int(M)$ converges to $y\in \P M$ and $\lim_{k\rightarrow \infty }|Du_k|=+\infty$. Since $M$ is $C^2$, then $gra(u)$ is a $C^2$ surface with $C^2$ boundary. As a result, by the continuity of normal vector fields, $gr(u)$ is tangent to $\P M$ at the point $(y,u(y))$. Without loss generality, we assume the normal vector of $gr(u)$ points outward of $M\PLH \R$ at $(y,u(y))$. By \eqref{eq:condition:mc}, it holds that $gr(\gamma)$ is equal to $\pm f(x)$ and $H_{\P M\PLH\R}\geq H_{gr(u)}$. From the Hopf Lemma, $\P M\PLH \R$ coincides with $\Sigma\PLH \R$ in a neighborhood of $(y, u(y))$. By the connectedness of $\P M$, $H_{\P M}=f$ on the whole $\P M$ or $H_{\P M}=-f$ on the whole $\P M$. This is a contradiction that $M$ is an {\nf} domain. Therefore $\sup_{x\in int(M)}\la Q, Q\ra<1$. \\
\indent The proof is complete. \\
\indent
\ep

\appendix
\section{Some results on mean curvature equations}
\indent In this section we derive the interior and boundary gradient estimates of the solution to a certain class of mean curvature equations.\\
\indent Throughout this section we suppose $\Omega$ is a $C^2$ bounded domain in a Riemannian manifold $N$ with dimension $n\geq 2$

Let $M$ be the product manifold $N\PLH\R$. Suppose $u$ is a function in $C^{3}(\Omega)$ satisfying $\omega=\sqrt{1+|Du|^2}$ and 
\be \label{eq:prescribed:mean:curvature}
-div(\F{Du}{\omega})+f(x,u,-\F{Du}{\omega}, \F{1}{\omega})=0
\ene 
where $f(x,z, X, r)$ is a $C^{1,\alpha}$  function on the tangent bundle of $M$ (written as $TM$), $p=(x,z)\in M$ for $x\in N,z\in  \R$, $v=X+r\P_r\in TM $ for $X\in TN$ and $r\in \R$.  Here $\P_r$ is the unit tangent vector filed of $\R$. 
\bt\label{thm:angle} Fix two positive numbers $c_0$ and $\mu$. Suppose (a) $u\in C^{3}(\Omega)$ satisfies \eqref{eq:prescribed:mean:curvature} and $\sup_{\Omega}|u|\leq c_0$;  (b) $\Sigma$ is the graph of $u$ over $\Omega$ and 
\be \label{constant:mu}
\begin{split}
\mathcal{N}_W(f)\leq \mu
\end{split}
\ene 
where $W=\{ (x,z):x\in \bar{\Omega}, |z|\leq c_0 \}
\subset M$ and $\mathcal{N}_{W}(f)$ is from Definition \ref{def:C:norm}.
Set $\Theta =\F{1}{\omega}$.  For any $\beta \in (0,1)$, there is a nonnegative constant $c_1$ only depending on $\mu$,  $\min_{x\in \bar{\Omega}}Ric_x$ and $\beta$  such that 
\be \label{formula:tde}
	\Delta \Theta +\beta |A|^2 \Theta+\F{\P f}{\P z}(1-\Theta^2)- c_1\Theta -\F{\P f}{\P  r}\la \nb \Theta, \P_r\ra\leq 0, 
	\ene 
where $A$ is the second fundamental form of $\Sigma$, $\Delta$ and $\nb$ is the Laplacian and the covariant derivative of $\Sigma$. 
\et
\bp By \cite[Lemma 2.2]{zhou19a} (its conclusion is true for any dimension greater than two) we have 
\be \label{de:A:hold}
\Delta \Theta+(|A|^2+\bar{R}ic(\vec{v},\vec{v}))\Theta-\la \nb H, \P_r\ra=0
\ene 
where $\vec{v}$ is the upward normal vector of $\Sigma$, $H$ is the mean curvature of $\Sigma$ {\wrt} $\vec{v}$ and $\bar{R}ic$ is the Ricci curvature of $M$. Notice that $\Theta =\la \vec{v},\P_r\ra$. Let $\P_r^T$ be the tangent component of $\P_r$ in $T\Sigma$. By definition we have 
\be 
\P_r^T=\P_r-\Theta\vec{v}
\ene 
 By  \eqref{eq:prescribed:mean:curvature}, $H=-f$ {\wrt} $\vec{v}$. As a result, we have  
\be\label{formula:theta}
\Delta \Theta+(|A|^2+\bar{R}ic(\vec{v},\vec{v})-\vec{v}(f))\Theta=0
\ene  
Let $\{\P_1, \cdots, \P_n \}$ be a local frame in $T\Omega$, $\sigma_{ij}=\la \P_i, \P_j\ra $, $(\sigma^{ij})=(\sigma_{ij})^{-1}$.  Let $u_{i}$ and $u_{ij}$ be the corresponding first and second covariant derivatives of $u$. Set $u^i=\sigma^{ik}u_k$. The gradient of $u$ is represented by $Du=u^k\P_k$. As a result,  $v^k=\F{u^k}{\omega}$,  $\vec{v}=\F{1}{\omega}\P_r-v^k\P_k$. Then $\{X_i =\P_i +u_i \P_r,i=1,\cdots, n\}$ is a frame on $T\Sigma$. The metric of $\Sigma$ is $(g_{ij})=(\la X_i, X_j\ra) =(\sigma_{ij}+u^iu^j)$ with its inverse martix $(g^{ij})=(\sigma^{ij}-\F{u^iu^j}{\omega^2})$. \\
\indent  The following two formulas are useful (see  \cite[section 2]{zhou18}). 
\begin{align*}
	v^k_i=\F{1}{\omega}g^{kl}u_{li}\quad |A|^2=v^k_iv_k^i 
\end{align*}
Let $X=p^1\P_1+\cdots+p^n\P_n$.  Let $f(x,z,X,r)$ be the function $f(x,z, p^1,\cdots, p^n, r)$.  We compute 
\be
\begin{split}
-\vec{v}(f) \Theta &=\Theta^2 \la D f, Du\ra =\Theta v^k \P_k f \\
   &=\Theta v^k(\P_k f+\F{\P f}{\P z} u_k+\sum_{j=1}^n\F{\P f}{\P p_{j}}v^j_k+\F{\P f}{\P r}\Theta_k)\\
   &=(v^k\P_k f+\F{\P f}{\P p_i}v^i_kv^k)\Theta +\Theta^2|Du|^2\F{\P f}{\P z}-\F{\P f}{\P  r}\la \nb\Theta, \P_r \ra \notag
\end{split}
\ene
In what follows we apply \eqref{constant:mu} (also see Remark \ref{rk:form}).  For any $\Sc>0$, we have $ 
\F{\P f}{\P p_i}v^i_kv^k\geq -\Sc |A|^2 -\F{\mu}{4\Sc}, \bar{R}ic(\vec{v},\vec{v})\geq Ric(\F{Du}{\omega},\F{Du}{\omega})\geq \min_{x\in \bar{\Omega}}Ric_x$. Notice that $\Theta^2|Du|^2=(1-\Theta^2)$. As a result, we have 
$$
(	\bar{R}ic(\vec{v}, \vec{v})-\vec{v}(f))\Theta \geq(\min_{x\in \bar{\Omega}}Ric_x -\mu(1+\F{1}{4\Sc}) -\Sc|A|^2)\Theta-\F{\P f}{\P r}\la \nb\Theta,\P_r\ra +\F{\P f}{\P z}(1-\Theta^2)
	$$
Putting this into \eqref{formula:theta} and letting $\Sc=1-\beta$ and $c_1=\max_{x\in \bar{\Omega}}|Ric_x| +\mu(1+\F{1}{4\Sc})$, we obtain the conclusion. 
\ep 
 Let $B_r(x)$ denote the embedded ball in $N$ centering at $x$ with radius $r$.  We roughly extend the result of Wang's estimate \cite[Theorem 1.1]{Wang98} into Riemannian manifolds as follows. 
\bt \label{thm:interior:estimate} Take the assumptions in Theorem \ref{thm:angle}. Suppose $\F{\P f}{\P z}\geq 0$.  \begin{enumerate}
	\item Let $B_{\rho}(x_0)\subset \Omega$. Then there is a positive constant $\mu_1$ depending on $\mu, c_0$ and the Ricci curvature on $B_{\rho}(x_0)$ such that   
\begin{gather}
	\max_{x\in B_{\F{1}{2}\rho}(x_0)}|Du|\leq \mu_1 \label{interior:estimate}
\end{gather}
\item Furthermore suppose $\Omega$ is $C^2$ bounded and $u\in C^{1}(\bar{\Omega})$. Then there is a positive constant $\mu_2$ depending on $\mu, c_0$ and the Ricci curvature on $\Omega$ such that 
\be \label{boundary:estimate:hold}
\max_{\bar{\Omega}}|Du|\leq  \mu_2(1+\max_{\P\Omega}|Du|)
\ene 
\end{enumerate}
\et 
\bp Recall that $\omega=\sqrt{1+|Du|^2}=\F{1}{\Theta}$. We just give the sketch of the proof. Since $\F{\P f}{\P z}\geq 0$, from \eqref{formula:tde} one sees that 
\be 
	\Delta \Theta - c_1\Theta -\F{\P f}{\P  r}\la \nb \Theta, \P_r\ra\leq 0, 
\ene 
 In the conclusion (1) set $q(x)=1+\F{u}{2c_0}-\F{3}{2\rho^2}dist^2(x_0, x)$, let $\eta=e^{Kq(x)}$ where $K$ is a sufficiently large constant only depending on $\rho,c_0$ and $\min_{x\in \bar{\Omega}}Ric_x$. Checking the maximum point of $\eta\omega$ will yield the estimate. \\
\indent As the conclusion (2), let $q(x)=u(x)$ and $\eta=e^{Ku}$ where $K$ is determined later.  If $\eta\omega$ achieves its maximum in $\Omega$ to obtain that 
\be 
e^{Ku}\omega \leq \max\{ C(K), e^{Kc_0}\sqrt{1+\max_{\P\Omega}|Du|^2}\}
\ene 
Here $K$ is a sufficiently large constant only depending on $\mu$ in \eqref{constant:mu}. This gives the conclusion (2). 
\ep 
Next, we record two boundary barrier functions used to construct the Perron subsolution and supsolution in section 4.3. Our derivation is similar as that of Marquardt \cite{Mar10}. 
\bt \label{lm:bet:function} Let $L(u)$ be the operator given by 
$$
L(u):=-div(\F{Du}{\omega})-H_1(x,u,-\F{Du}{\omega}, \F{1}{\omega})+H_{2}(x,u,-\F{Du}{\omega} , \F{1}{\omega})\F{1}{\omega}$$
where $H_i(x, z, X, r)$ for $i=1,2$ are $C^{1,\alpha}$ functions in the tangent bundle of $N\PLH\R$. Let $d(x)=\inf_{y\in\P\Omega}d(x,y)$.  Fix $\vp \in C^2(\bar{\Omega})$ and any positive constant $c_0>\max_{\bar{\Omega}}|\vp|$.  
\begin{enumerate}
\item Suppose 
\be \label{condition:restriction:A}
H_{\P\Omega}(x)\geq  H_1(x,\vp,  \gamma(x), 0)
\ene
where $\gamma(x)$ is the outward normal vector of $\P\Omega$. Then there are two positive constants $\kappa>0$, $\nu$  depending on $c_0$, $H_1, H_2$ and the $C^2$ norm of $\vp$ and $d$ such that for $u_{+}(x)=\vp + \F{\log(1+\kappa d(x))}{\nu}$ it holds that 
\be 
   L(u_{+})\geq 0 
\ene
on the domain $\Gamma_{+}:=\{ x\in \Omega:  u_{+}(x)<c_0\}$ with $\pm u_{\pm}=c_0$ on $\P\Gamma_{+}\backslash \P\Omega$;
\item  Suppose 
\be \label{condition:restriction:B}
H_{\P\Omega}(x)\geq  -H_1(x,\vp,  -\gamma(x), 0)
\ene
then Then there are two positive constants $\kappa>0$, $\nu$  depending on $c_0$, $H_1, H_2$ and the $C^2$ norm of $\vp$ and $d$ such that for $u_{+}(x)=\vp + \F{\log(1+\kappa d(x))}{\nu}$ it holds that 
\be 
-L(u_{-})\geq 0 
\ene
on the domain $\Gamma_{-}:=\{ x\in \Omega:  u_{-}(x)<c_0\}$ with $- u_{-}=c_0$ on $\P\Gamma_{-}\backslash \P\Omega$;
\end{enumerate}   
\et
\br To construct the subsolution and supersolution of $L(u)=0$ on $\Omega$, we do not require that \eqref{condition:restriction:A}-\eqref{condition:restriction:B} strictly holds on $\P \Omega$.  
\er  
\bp Let $\{e_1,\cdots, e_n\}$ be any local orthonormal frame of $T N$. As a result any  vector field  $p$ in $TN$ can be rewritten as $p=p_ie_i$. Thus $|p|^2=p_i^2$. Then we follow the notation in \cite[(14.3)]{GT01}. For any $u\in C^2(\bar{\Omega})$ we write \be 
L'(u):=-(1+|Du|^2)^{\F{3}{2}} L(u)=a^{ij}u_{ij}+b
\ene 
where the arguments of all functions are $x,u ,p = Du$, $\omega=\sqrt{1+|p|^2}$,  $\Lambda=1+|p|^2$ and 
\be 
\begin{split}
a^{ij}(x,u,p)&=\Lambda a_{\infty}^{ij}+a_0^{ij},\quad b(x,u,p)=|p|\Lambda b_\infty+b_0\\
a_\infty^{ij}&=\delta_{ij}-\F{p_ip_j}{|p|^2},\quad a_0^{ij}=\F{p_ip_j}{|p|^2}, b_\infty=-H_1(x,u, -\F{p}{\omega},\F{1}{\omega}),\\  
b_0(x,u,p)&= -\Lambda (\F{1}{\omega+|p|} H_1(x,u, -\F{p}{\omega},\F{1}{\omega})+H_2(x,u, -\F{p}{\omega},\F{1}{\omega}))
\end{split}
\ene 
 Now define $\phi(r)=\F{\log(1+\kappa r)}{\nu}$. It holds that 
\be
\phi'(r)=\F{\kappa}{\nu(1+\kappa r )},\quad  \F{\phi''}{(\phi')^2}=-\nu
\ene 
Let $d(x)$ be the function $\inf_{y\in\P\Omega}d(x,y)$ for $x$ in $\Omega$. We define $u_{\pm}(x)=\vp(x) \pm \phi(d(x))$. Define  $\Gamma:=\{x\in \Omega: d(x)<d_0\}$ where $d_0$ is a sufficiently small positive constant determined later such that for any $x$ in $\Gamma_{\pm}$ there is a unique point $y\in \P \Omega$ such that $d(x,y)=d(x,\P\Omega)$. Let $\gamma(y)$ be the outward normal vector of $y$ which is equal to $-Dd(y)$. From now on let $u=u_{\pm}$ and $p$ be $Du_{\pm}$.\\
\indent  In what follows we use $q(s)$ denote different  continuous functions satisfying $q(0)=0$, $|\F{q(r)}{r}|\leq C$ for all $r\in [-d_0,d_0]$ and some positive constant $C$. Let $\mu$ be the constant 
\be 
\max_{x\in \bar{\Omega}, |z|\leq c_0, \la X, X\ra\leq 1, |r|\leq 1}\{|H_1(x,z,X,r)|+|H_2(x,z, X,r)|\}
\ene 
Let $\Delta$ be the Laplacian on $N$. Then on $\Gamma$ it holds 
\be 
\begin{split}
  |b_0(x,u_{\pm},Du_{\pm})|&\leq \mu \Lambda, \quad  \Delta d(x)=-H_{\P\Omega}(y)+q(d(x))\\
  H_{1}(x, u_{\pm}(x), -\F{p}{\omega}, \F{1}{\omega})&=H_{1}(y, \vp(y), \pm \gamma(y), 0)+q(d(x)+|\phi'|^{-2}(d(x)))\\
  p&=D\vp\pm \phi'(d(x))Dd=\phi'(d(x))(1+q(\F{1}{\phi'(d(x))})\\
  a^{ij}(x,u_{\pm},p)(\phi')^2d_id_j &=\Lambda (1+q(\F{1}{\phi'(d(x))}))
\end{split}
\ene 
 As a result 
\be	\label{key:step}
\begin{split}
\pm L'(u)(x)&=\Lambda \phi'\Delta d(x)\pm a^{ij}\vp_{ij}
+\phi'(1+q(\F{1}{\phi'}))\Lambda \{\mp H_{1}(y, \vp(y),\pm\gamma(y), 0)\\
&+q(d+|\phi'|^{-2})\}+b_0-\nu \Lambda(1+q(\F{1}{\phi'(d)}))\\
&\leq \Lambda \phi'(-H_{\P\Omega}(y)\mp H_{1}(y,\vp(y), \pm \gamma(y), 0))\\
&+\Lambda(C\phi'd+q(\F{1}{\phi'})- \nu)
\end{split} 
\ene 
Since $\nu\geq 1$, we have $\phi'd=\F{1}{\nu}\F{\kappa d}{1+\kappa d}\leq 1$. Now there are two positive constants $C_2\geq 1,C_3$ depending $\mu, d(x)$ and $\vp$ such that  
if 
\be \label{key:two:constant:restriction}
\phi'\geq C_2, \nu\geq  \max\{1, C_3\}
\ene 
then for any $x$ in $ \Gamma$ we have 
\be \label{dq:det}
\pm L'(u)(x)\leq \Lambda \phi'(-H_{\P\Omega}(y)\mp H_{1}(y,\vp(y),\pm \gamma(y), 0))\leq 0
\ene 
In the last inequality, we use the condition \eqref{condition:restriction:A} and \eqref{condition:restriction:B}. Thus $\pm Lu_{\pm}\geq 0$ on $\Gamma$. \\
\indent Now to achieve the goal in \eqref{dq:det}, first fix $\nu\geq \max\{1, C_3\}$. For $x$ in $ \Gamma$, we require that 
\be 
\begin{split}
\phi'(d(x))&=\F{\kappa }{\nu(1+\kappa d(x)) }\geq\F{\kappa }{\nu(1+\kappa d_0) }\geq C_2;\\
\phi(d_0)&=\F{\log(1+ \kappa d_0)}{\nu}\geq c_0+\max_{\bar{\Omega}}|\vp| :=C_4
\end{split}
\ene 
The above two inequalities are satisfied by choosing $d_0$ small such that $d_0C_2\nu \leq \F{1}{2}$ and choosing $\kappa$ such that 
\be 
\kappa =\max\{ \F{C_2\nu}{1-C_2\nu d_0}, \F{1}{d_0}(e^{C_4\nu}-1)\}
\ene 
\indent With those fixed $\kappa, \nu$ one sees that $\pm L(u_\pm)\geq 0$ on $\Gamma=\{x:d(x)<d_0\}$. Now define $\Gamma_{\pm}=\{x\in \Gamma: \pm u_{\pm}(x)<c_0\}$. We obtain the desirable conclusion. 
\ep 

  \bibliographystyle{plain}
  \bibliography{Ref-MC}
 \end{document}